\documentclass[12pt]{preprint} 

\usepackage[a4paper,innermargin=1.2in,outermargin=1.2in,
bottom=1.5in,marginparwidth=1in,marginparsep
=3mm]{geometry}

\usepackage[T1]{fontenc} % Use 8-bit encoding that has 256 glyphs
\usepackage{fourier} % Use the Adobe Utopia font for the document - comment this line to return to the LaTeX default
\usepackage[english]{babel} % English language/hyphenation
\usepackage{amsmath,amsfonts,amsthm} % Math packages
\usepackage{tikz-cd}
\usepackage{lipsum} % Used for inserting dummy 'Lorem ipsum' text into the template
\usepackage{shuffle}
\usepackage{ amssymb }
\usepackage{fancyhdr} % Custom headers and footers
\usepackage{abstract}
\usepackage{comment}
\usepackage{stmaryrd}
\usepackage{cprotect}
\usepackage{hyperref}
\usepackage{mathrsfs}
\usepackage{mhequ}
\usepackage{orcidlink}

\colorlet{darkblue}{blue!90!black}
\colorlet{darkred}{red!90!black}

\newcommand{\fraks}{\mathfrak{s}}

\newcommand{\cC}{\mathcal{C}}

\newcommand{\RR}{\mathbb{R}}
\newcommand{\R}{\mathbb{R}}

\newcommand{\NN}{\mathbb{N}}

\newcommand{\bfPi}{\pmb{\Pi}}

\newcommand{\bfK}{\pmb{\mathcal{K}}}

\newcommand{\opP}{\operatorname{P}}

\renewcommand{\bf}{\mathbf{f}}

\newcommand{\eps}{\varepsilon}

\newcommand{\E}{\mathbb{E}}

\newcommand{\vertiii}[1]{{\left\vert\kern-0.25ex\left\vert\kern-0.25ex\left\vert #1 
		\right\vert\kern-0.25ex\right\vert\kern-0.25ex\right\vert}}

\newtheorem{theorem}{Theorem}[section]
\newtheorem{corollary}{Corollary}[theorem]

\newtheorem{lemma}[theorem]{Lemma}
\newtheorem{definition}{Definition}
\newtheorem{prop}[theorem]{Proposition}
\newtheorem{assumption}[theorem]{Assumption}
\theoremstyle{remark}
\newtheorem{remark}[theorem]{Remark}

\numberwithin{equation}{section} % Number equations within sections (i.e. 1.1, 1.2, 2.1, 2.2 instead of 1, 2, 3, 4)
\numberwithin{figure}{section} % Number figures within sections (i.e. 1.1, 1.2, 2.1, 2.2 instead of 1, 2, 3, 4)
\numberwithin{table}{section} % Number tables within sections (i.e. 1.1, 1.2, 2.1, 2.2 instead of 1, 2, 3, 4)

\colorlet{symbols}{blue!90!black}
\colorlet{testcolor}{green!60!black}

\usetikzlibrary{calc}
\usetikzlibrary{shapes.misc}
\usetikzlibrary{shapes.symbols}
\usetikzlibrary{shapes.geometric}
\usetikzlibrary{snakes}
\usetikzlibrary{decorations}
\usetikzlibrary{decorations.markings}

\tikzset{
	root/.style={circle,fill=testcolor,inner sep=0pt, minimum size=2mm},
	broot/.style={circle,fill=gray,inner sep=0pt, minimum size=2mm},
	dot/.style={circle,fill=black,inner sep=0pt, minimum size=1mm},
		reddot/.style={circle,fill=red,inner sep=0pt, minimum size=1mm},
			bluedot/.style={circle,fill=blue,inner sep=0pt, minimum size=1mm},
	eps/.style={circle,fill=white,draw=symbols,inner sep=0pt,minimum size=0.8mm},
	int/.style={circle,fill=black,draw=black,inner sep=0pt,minimum size=0.7mm},
	var/.style={circle,fill=black!10,draw=black,inner sep=0pt, minimum size=2mm},
	dotred/.style={circle,fill=black!50,inner sep=0pt, minimum size=2mm},
	generic/.style={semithick,shorten >=1pt,shorten <=1pt},
	dist/.style={ultra thick,draw=testcolor,shorten >=1pt,shorten <=1pt},
	testfcn/.style={ultra thick,testcolor,shorten >=1pt,shorten <=1pt,<-},
	testfcnx/.style={ultra thick,testcolor,shorten >=1pt,shorten <=1pt,<-,
		postaction={decorate,decoration={markings,mark=at position 0.6 with {\drawx}}}},
	keps/.style={semithick,shorten >=1pt,shorten <=1pt,densely dashed,->},
	kprimex/.style={semithick,shorten >=1pt,shorten <=1pt,densely dashed,->,
		postaction={decorate,decoration={markings,mark=at position 0.4 with {\drawx}}}},
	kernel/.style={semithick,shorten >=1pt,shorten <=1pt,->},
	multx/.style={shorten >=1pt,shorten <=1pt,
		postaction={decorate,decoration={markings,mark=at position 0.5 with {\drawx}}}},
	kernelx/.style={semithick,shorten >=1pt,shorten <=1pt,->,
		postaction={decorate,decoration={markings,mark=at position 0.4 with {\drawx}}}},
	kernel1/.style={->,semithick,shorten >=1pt,shorten <=1pt,postaction={decorate,decoration={markings,mark=at position 0.45 with {\draw[-] (0,-0.1) -- (0,0.1);}}}},
	kernel2/.style={->,semithick,shorten >=1pt,shorten <=1pt,postaction={decorate,decoration={markings,mark=at position 0.45 with {\draw[-] (0.05,-0.1) -- (0.05,0.1);\draw[-] (-0.05,-0.1) -- (-0.05,0.1);}}}},
	kernelBig/.style={semithick,shorten >=1pt,shorten <=1pt,decorate, decoration={zigzag,amplitude=1.5pt,segment length = 3pt,pre length=2pt,post length=2pt}},
	rho/.style={dotted,semithick,shorten >=1pt,shorten <=1pt},
	renorm/.style={shape=circle,fill=white,inner sep=1pt},
	labl/.style={shape=rectangle,fill=white,inner sep=1pt},
	xi/.style={circle,fill=symbols!10,draw=symbols,inner sep=0pt,minimum size=1.2mm},
	xix/.style={crosscircle,fill=symbols!10,draw=symbols,inner sep=0pt,minimum size=1.2mm},
	xib/.style={circle,fill=symbols!10,draw=symbols,inner sep=0pt,minimum size=1.6mm},
	xibx/.style={crosscircle,fill=symbols!10,draw=symbols,inner sep=0pt,minimum size=1.6mm},
	not/.style={circle,fill=symbols,draw=symbols,inner sep=0pt,minimum size=0.5mm},
cumu2n/.style={inner sep=3pt},
cumu2/.style={draw=red!80,fill=red!40},
cumu2b/.style={draw=blue!80,fill=blue!40},
cumu2nv/.style={inner sep=3pt},
cumu2v/.style={draw=red!80,fill=white,very thick},
cumu3/.style={regular polygon, regular polygon sides=3,draw=red!80,rounded corners=3pt,fill=red!40,minimum size=5mm},
cumu4/.style={regular polygon, regular polygon sides=4,draw=red!80,rounded corners=3pt,fill=red!40,minimum size=7mm},
cumu5/.style={regular polygon, regular polygon sides=5,draw=red!80,rounded corners=3pt,fill=red!40,minimum size=7mm},
	>=stealth,
	not/.style={circle,fill=symbols,draw=symbols,inner sep=0pt,minimum size=0.5mm},
kernels2/.style={very thick,segment length=12pt},
arrho/.style={dotted,semithick,shorten >=1pt,shorten <=1pt, ->},
	}
\makeatletter

\def\DeclareSymbol#1#2#3{\expandafter\gdef\csname MH@symb@#1\endcsname{\tikz[baseline=#2,scale=0.15,draw=symbols]{#3}}}
\def\<#1>{\csname MH@symb@#1\endcsname}
\makeatother
\def\DeclareSymbol#1#2#3{%
	\expandafter\gdef\csname MH@symb@#1\endcsname{\tikzsetnextfilename{symbol#1}%
		\tikz[baseline=#2,scale=0.15,draw=symbols,line join=round]{#3}}%
	\expandafter\gdef\csname MH@symb@#1s\endcsname{\scalebox{0.75}{\tikzsetnextfilename{symbol#1}%
			\tikz[baseline=#2,scale=0.15,draw=symbols,line join=round]{#3}}}%
	\expandafter\gdef\csname MH@symb@#1ss\endcsname{\scalebox{0.65}{\tikzsetnextfilename{symbol#1}%
			\tikz[baseline=#2,scale=0.15,draw=symbols,line join=round]{#3}}}%
}

\makeatletter

\def\DeclareSymbol#1#2#3{\expandafter\gdef\csname MH@symb@#1\endcsname{\tikz[baseline=#2,scale=0.15,draw=symbols]{#3}}}
\def\<#1>{\csname MH@symb@#1\endcsname}
\makeatother
\DeclareSymbol{X}{-2.4}{\node[eps] {};}
\DeclareSymbol{I}{0}{\draw[white] (-.5,0) -- (.5,0); \draw (0,0)  -- (0,1.5) {};}
\DeclareSymbol{bI}{0}{
\draw[white] (-.5,0) -- (.5,0); \draw[kernels2] (0,0)  -- (0,1.5);
}

\DeclareSymbol{1}{0}{\draw[white] (-.5,0) -- (.5,0); \draw (0,0)  -- (0,1.5) node[eps] {};}
\DeclareSymbol{2}{0}{\draw (-0.5,1.5) node[eps] {} -- (0,0) -- (0.5,1.5) node[eps] {};}
\DeclareSymbol{3}{0}{\draw (0,0) -- (0,1.5) node[eps] {}; \draw (-.7,1.3) node[eps] {} -- (0,0) -- (.7,1.3) node[eps] {};}
\DeclareSymbol{E3}{0}{\draw (0,0) -- (0,1.5) node[not] {}; \draw (-.7,1.3) node[not] {} -- (0,0) -- (.7,1.3) node[not] {}; \node[eps] (0,0) {};}

\DeclareSymbol{E5E4}{1}{
\draw (.5,2) node[not] {} -- (1.5,1.5) -- (2.5,2) node[not] {};
\draw (1.5,1.5) -- (1.5,2.7) node[not] {}; \draw (0.8,2.5) node[not] {} -- (1.5,1.5) -- (2.2,2.5) node[not] {};
\draw (0,0) -- (0,1.5) node[not] {}; \draw (-.7,1.3) node[not] {} -- (0,0) -- (.7,1.3) node[not] {};
\draw (-1.2,0.7) node[not] {} -- (0,0);
\draw (0,0) to[bend right=50] (1.5,1.5); 
\node[eps] at (0,0) {};\node[eps] at (1.5,1.5) {};
}
\DeclareSymbol{E4E4}{1}{
\draw (.5,2) node[not] {} -- (1.5,1.5) -- (2.5,2) node[not] {};
\draw (1,2.6) node[not] {} -- (1.5,1.5) -- (2,2.6) node[not] {};
\draw (0,0) -- (0,1.5) node[not] {}; \draw (-.7,1.3) node[not] {} -- (0,0) -- (.7,1.3) node[not] {};
\draw (-1.2,0.7) node[not] {} -- (0,0);
\draw (0,0) to[bend right=50] (1.5,1.5); 
\node[eps] at (0,0) {};\node[eps] at (1.5,1.5) {};
}
\DeclareSymbol{3E4}{1}{
\draw (1.5,1.5) -- (1.5,2.7) node[not] {}; \draw (0.8,2.5) node[not] {} -- (1.5,1.5) -- (2.2,2.5) node[not] {};
\draw (0,0) -- (0,1.5) node[not] {}; \draw (-.7,1.3) node[not] {} -- (0,0) -- (.7,1.3) node[not] {};
\draw (-1.2,0.7) node[not] {} -- (0,0);
\draw (0,0) to[bend right=50] (1.5,1.5); 
\node[eps] at (0,0) {};
}

\DeclareSymbol{3E3}{1}{
\draw (1.5,1.5) -- (1.5,2.7) node[not] {}; \draw (0.8,2.5) node[not] {} -- (1.5,1.5) -- (2.2,2.5) node[not] {};
\draw (0,0) -- (0,1.5) node[not] {}; \draw (-.7,1.3) node[not] {} -- (0,0) -- (.7,1.3) node[not] {};
\draw (0,0) to[bend right=50] (1.5,1.5); 
\node[eps] at (0,0) {};
}

\DeclareSymbol{E5}{0}{\draw (0,0) -- (0,1.5) node[not] {}; \draw (-.7,1.3) node[not] {} -- (0,0) -- (.7,1.3) node[not] {};
\draw (-1.2,0.7) node[not] {} -- (0,0) -- (1.2,0.7) node[not] {}; \node[eps] (0,0) {};}
\DeclareSymbol{E52}{-3}{\draw (0,0) -- (0,-1); \draw (1.5,-0.4) node[not] {} -- (0,-1) -- (-1.5,-0.4) node[not] {}; 
\draw (-1,0.5) node[not] {} -- (0,0) -- (1,0.5) node[not] {};
\draw (0,0) -- (0,1.2) node[not] {}; \draw (-.7,1) node[not] {} -- (0,0) -- (.7,1) node[not] {};\node[eps] (0,0) {};}
\DeclareSymbol{E32}{-3}{\draw (0,0) -- (0,-1); \draw (1.5,-0.4) node[not] {} -- (0,-1) -- (-1.5,-0.4) node[not] {}; 
\draw (0,0) -- (0,1.2) node[not] {}; \draw (-.7,1) node[not] {} -- (0,0) -- (.7,1) node[not] {};\node[eps] (0,0) {};}
\DeclareSymbol{E50}{-3}{\draw (0,0) -- (0,-1); 
\draw (-1,0.5) node[not] {} -- (0,0) -- (1,0.5) node[not] {};
\draw (0,0) -- (0,1.2) node[not] {}; \draw (-.7,1) node[not] {} -- (0,0) -- (.7,1) node[not] {};\node[eps] (0,0) {};}

\DeclareSymbol{3E2}{-3}{\draw (0,0) -- (0,-1); \draw (1,0) node[not] {} -- (0,-1) -- (-1,0) node[not] {}; \draw (0,0) -- (0,1.2) node[not] {}; \draw (-.7,1) node[not] {} -- (0,0) -- (.7,1) node[not] {};\node[eps] at (0,-1) {};}

\DeclareSymbol{31}{-3}{\draw (0,0) -- (0,-1) -- (1,0) node[eps] {}; \draw (0,0) -- (0,1.2) node[eps] {}; \draw (-.7,1) node[eps] {} -- (0,0) -- (.7,1) node[eps] {};}
\DeclareSymbol{30}{-3}{\draw (0,0) -- (0,-1); \draw (0,0) -- (0,1.2) node[eps] {}; \draw (-.7,1) node[eps] {} -- (0,0) -- (.7,1) node[eps] {};}
\DeclareSymbol{32}{-3}{\draw (0,0) -- (0,-1); \draw (1,0) node[eps] {} -- (0,-1) -- (-1,0) node[eps] {}; \draw (0,0) -- (0,1.2) node[eps] {}; \draw (-.7,1) node[eps] {} -- (0,0) -- (.7,1) node[eps] {};}
\DeclareSymbol{32b}{-4.3}{\draw (0,0.5) -- (0,-1); \draw (1,0) node[not] {} -- (0,-1) -- (-1,0) node[not] {};}
\DeclareSymbol{22}{-3}{\draw (0,0.3) -- (0,-1); \draw (1,0) node[eps] {} -- (0,-1) -- (-1,0) node[eps] {};\draw (-.7,1) node[eps] {} -- (0,0.3) -- (.7,1) node[eps] {};}
\DeclareSymbol{20}{-3}{\draw (0,0) -- (0,-1);\draw (-.7,1) node[not] {} -- (0,0) -- (.7,1) node[not] {};}
\DeclareSymbol{12}{-3}{\draw (0,0.3) -- (0,-1); \draw (1,0) node[not] {} -- (0,-1) -- (-1,0) node[not] {};\draw (-.7,1) node[not] {} -- (0,0.3);}
\DeclareSymbol{10}{-3}{\draw (0,0.3) -- (0,-1);\draw (-.7,1) node[not] {} -- (0,0.3);}
\DeclareSymbol{21}{-3}{\draw (0,0.3) -- (0,-1) -- (1,0) node[not] {};\draw (-.7,1) node[not] {} -- (0,0.3) -- (.7,1) node[not] {};}

\DeclareSymbol{1'}{0}{\draw[white] (-.5,0) -- (.5,0); \draw (0,0)  -- (0,1.5) node[not] {}; \draw ( -0.2,1)-- ( 0.2, 1);}

\DeclareSymbol{3'}{0}{\draw (0,0) -- (0,1.5) node[not] {};
\draw  (-0.25,1.1) -- ( 0.25, 1.1);
 \draw (-.7,1.3) node[eps] {} -- (0,0) -- (.7,1.3) node[eps] {};}

\DeclareSymbol{3''}{0}{\draw (0,0) -- (0,1.5) node[eps] {};
 \draw (-.7,1.3) node[not] {} -- (0,0) -- (.7,1.3) node[not] {};
 \draw  (-0.8,0.7) -- ( -0.25, 1);
  \draw  (0.8,0.7) -- ( 0.25, 1);
 }
 
 \DeclareSymbol{3'''}{0}{\draw (0,0) -- (0,1.5) node[not] {};
\draw  (-0.3,1.05) -- ( 0.3, 1.05);
 \draw (-.7,1.3) node[not] {} -- (0,0) -- (.7,1.3) node[not] {};
 \draw  (-0.8,0.7) -- ( -0.25, 1);
  \draw  (0.8,0.7) -- ( 0.25, 1);
 }
 
 \DeclareSymbol{3'2}{-3}{\draw (0,0) -- (0,-1); \draw (1,0) node[eps] {} -- (0,-1) -- (-1,0) node[eps] {}; \draw (0,0) -- (0,1.4) node[not] {}; \draw (-.7,1) node[eps] {} -- (0,0) -- (.7,1) node[eps] {};
\draw  (-0.25,1) -- ( 0.25, 1); 
 }
 
 \DeclareSymbol{3''2}{-3}{\draw (0,0) -- (0,-1); \draw (1,0) node[eps] {} -- (0,-1) -- (-1,0) node[eps] {}; \draw (0,0) -- (0,1.3) node[eps] {}; \draw (-.7,1) node[not] {} -- (0,0) -- (.7,1) node[not] {};
 \draw  (-0.8,0.5) -- ( -0.25, 0.8);
  \draw  (0.8,0.5) -- ( 0.25, 0.8); 
 }
 
  \DeclareSymbol{3'''2}{-3}{\draw (0,0) -- (0,-1); \draw (1,0) node[eps] {} -- (0,-1) -- (-1,0) node[eps] {}; \draw (0,0) -- (0,1.4) node[not] {}; \draw (-.7,1) node[not] {} -- (0,0) -- (.7,1) node[not] {};
 \draw  (-0.8,0.5) -- ( -0.25, 0.8);
  \draw  (0.8,0.5) -- ( 0.25, 0.8); 
  \draw ( -0.25,1)-- ( 0.25, 1);
 }
 
  \DeclareSymbol{3'2'}{-3}{\draw (0,0) -- (0,-1); \draw (1,0) node[eps] {} -- (0,-1) -- (-1,0) node[not] {}; \draw (0,0) -- (0,1.4) node[not] {}; \draw (-.7,1) node[eps] {} -- (0,0) -- (.7,1) node[eps] {};
\draw  (-0.25,1) -- ( 0.25, 1); 
\draw  (-0.8,-0.6) -- ( -0.4, -0.2);
 }
 
  \DeclareSymbol{3'2''}{-3}{\draw (0,0) -- (0,-1); \draw (1,0) node[not] {} -- (0,-1) -- (-1,0) node[not] {}; \draw (0,0) -- (0,1.4) node[not] {}; \draw (-.7,1) node[eps] {} -- (0,0) -- (.7,1) node[eps] {};
\draw  (-0.25,1) -- ( 0.25, 1); 
\draw  (-0.8,-0.6) -- ( -0.4, -0.2);
\draw  (0.8,-0.6) -- ( 0.4, -0.2);
 }

   \DeclareSymbol{3''2'}{-3}{\draw (0,0) -- (0,-1); \draw (1,0) node[eps] {} -- (0,-1) -- (-1,0) node[not] {}; \draw (0,0) -- (0,1.3) node[eps] {}; \draw (-.7,1) node[not] {} -- (0,0) -- (.7,1) node[not] {};
 \draw  (-0.8,0.5) -- ( -0.25, 0.8);
  \draw  (0.8,0.5) -- ( 0.25, 0.8); 
\draw  (-0.8,-0.6) -- ( -0.4, -0.2);
 }
  \DeclareSymbol{32'}{-3}{\draw (0,0) -- (0,-1); \draw (1,0) node[eps] {} -- (0,-1) -- (-1,0) node[not] {}; \draw (0,0) -- (0,1.3) node[eps] {}; \draw (-.7,1) node[eps] {} -- (0,0) -- (.7,1) node[eps] {};
\draw  (-0.8,-0.6) -- ( -0.4, -0.2);
 }
  \DeclareSymbol{3'0}{-3}{\draw (0,0) -- (0,-1); 
   \draw (0,0) -- (0,1.4) node[not] {}; \draw (-.7,1) node[eps] {} -- (0,0) -- (.7,1) node[eps] {};
\draw  (-0.25,1) -- ( 0.25, 1); 
 }

\DeclareSymbol{1}{0}{\draw[white] (-.4,0) -- (.4,0); \draw (0,0)  -- (0,1.2) node[eps] {};}
\DeclareSymbol{2}{0}{\draw (-0.5,1.2) node[eps] {} -- (0,0) -- (0.5,1.2) node[eps] {};}
\DeclareSymbol{11}{0}{\draw (0,1.8) node[eps] {} -- (-0.7,0.9) -- (0,0) -- (0.7,1) node[eps] {};}
\DeclareSymbol{10}{0}{\draw (0,1.8) node[eps] {} -- (-0.8,0.9) -- (0,0);}
\DeclareSymbol{21}{0.7}{\draw (-1,1.8) node[eps] {} -- (-0.5,0.9); \draw (0,1.8) node[eps] {} -- (-0.5,0.9) -- (0,0) -- (0.5,0.9) node[eps] {};}
\DeclareSymbol{20}{0.7}{\draw (-1,1.8) node[eps] {} -- (-0.5,0.9); \draw (0,1.8) node[eps] {} -- (-0.5,0.9) -- (-0.5,0);}
\DeclareSymbol{210}{1.1}{\draw (-0.5,2.4) node[eps] {} -- (-1,1.6);\draw (-1.5,2.4) node[eps] {} -- (-0.5,0.8); \draw (0,1.6) node[eps] {} -- (-0.5,0.8) -- (-0.5,0);}
\DeclareSymbol{211}{1.1}{\draw (-0.5,2.4) node[eps] {} -- (-1,1.6);\draw (-1.5,2.4) node[eps] {} -- (-0.5,0.8); \draw (0,1.6) node[eps] {} -- (-0.5,0.8) -- (0,0) -- (0.5,0.8) node[eps] {};}
\DeclareSymbol{22j}{0.7}{\draw (-1.5,1.8) node[eps] {} -- (-1,0.9) -- (0,0) -- (1,0.9) -- (1.5,1.8) node[eps] {};
\draw (-0.5,1.8) node[eps] {} -- (-1,0.9);\draw (0.5,1.8) node[eps] {} -- (1,0.9);}

\DeclareSymbol{2'}{0}{\draw (-0.65,1.3) node[not] {} -- (0,0) -- (0.6,1.3) node[eps] {};
\draw  (-0.8,0.7) -- ( -0.15, 1.1);} 

\DeclareSymbol{2''}{0}{\draw (-0.65,1.3) node[not] {} -- (0,0) -- (0.6,1.3) node[not] {};
\draw  (-0.8,0.7) -- ( -0.15, 1.1);
\draw  (0.8,0.7) -- ( 0.15, 1.1);
} 
 
\DeclareSymbol{2'1}{0.7}{\draw (-1,1.8) node[not] {} -- (-0.5,0.9); \draw (0,1.8) node[eps] {} -- (-0.5,0.9) -- (0,0) -- (0.5,0.9) node[eps] {};
\draw  (-1.1,1.2) -- ( -0.5, 1.6);} 

\DeclareSymbol{2''1}{0.7}{\draw (-1,1.8) node[not] {} -- (-0.5,0.9); \draw (0,1.8) node[not] {} -- (-0.5,0.9) -- (0,0) -- (0.5,0.9) node[eps] {};
\draw  (-1.1,1.2) -- ( -0.5, 1.6);
\draw  (.15,1.2) -- ( -0.45, 1.6);} 
 
\DeclareSymbol{2''1'}{0.7}{\draw (-1,1.8) node[not] {} -- (-0.5,0.9); \draw (0,1.8) node[not] {} -- (-0.5,0.9) -- (0,0) -- (0.55,1) node[not] {};
\draw  (-1.1,1.2) -- ( -0.5, 1.6);
\draw  (.15,1.2) -- ( -0.45, 1.6);
\draw  (0.8,0.4) -- ( 0, 0.8);
} 
\DeclareSymbol{2'1'}{0.7}{\draw (-1,1.8) node[not] {} -- (-0.5,0.9); \draw (0,1.8) node[eps] {} -- (-0.5,0.9) -- (0,0) -- (0.55,1) node[not] {};
\draw  (-1.1,1.2) -- ( -0.5, 1.6);
\draw  (0.8,0.4) -- ( 0, 0.8);} 

\DeclareSymbol{21'}{0.7}{\draw (-1,1.8) node[eps] {} -- (-0.5,0.9); \draw (0,1.8) node[eps] {} -- (-0.5,0.9) -- (0,0) -- (0.55,1) node[not] {};
\draw  (0.8,0.4) -- ( 0, 0.8);} 

\DeclareSymbol{2'0}{0.7}{\draw (-1,1.8) node[not] {} -- (-0.5,0.9); \draw (0,1.8) node[eps] {} -- (-0.5,0.9) -- (-0.5,0);
\draw  (-1.1,1.2) -- ( -0.5, 1.6);}

\DeclareSymbol{b1}{0}{
\draw[white] (-.5,0) -- (.5,0); \draw[kernels2] (0,0)  -- (0,1.5);
 \draw  (0,1.5) node[eps] {};
}

\DeclareSymbol{b2}{0}{
\draw[kernels2] (0,0) -- (-0.5,1.4);
\draw[kernels2] (0,0) -- (+0.5,1.4); 
 \draw  (-0.5,1.4) node[eps] {};
 \draw  (0.5,1.4) node[eps] {};
 }
 
 \DeclareSymbol{20b}{0.7}{\draw (-1,1.8) node[eps] {} -- (-0.5,0.9);
  \draw (0,1.8) node[eps] {} -- (-0.5,0.9);
  \draw[kernels2](-0.5,0.9) -- (-0.5,0);}
  
\DeclareSymbol{210b}{1.1}{\draw (-0.5,2.4) node[eps] {} -- (-1,1.6);\draw (-1.5,2.4) node[eps] {} -- (-0.5,0.8); \draw (0,1.6) node[eps] {} -- (-0.5,0.8);
\draw[kernels2] (-0.5,0.8) -- (-0.5,0);}
 
\DeclareSymbol{bI}{0}{\draw[white] (-.5,0) -- (.5,0); \draw[kernels2] (0,0)  -- (0,1.5) {};} 

\DeclareSymbol{Xi2}{0}{\draw[white] (-.5,0) -- (.5,0); \draw (0,0)  -- (0,1.5) node[eps] {};
\draw (0,0) node[eps] {};}

 \DeclareSymbol{1'Xi}{0}{\draw[white] (-.5,0) -- (.5,0); \draw (0,0)  -- (0,1.5) node[not] {}; \draw ( -0.2,1)-- ( 0.2, 1);
 \draw (0,0) node[eps] {};
 }
\def\CCE{\mathbb{E}}
\def\CCV{\mathbb{V}}
\def\CCG{\mathbb{G}}

\newcommand{\horrule}[1]{\rule{\linewidth}{#1}} % Create horizontal rule command with 1 argument of height

\title{	
\normalfont \normalsize 
\horrule{0.5pt} \\[0.4cm] % Thin top horizontal rule
\large Canonical solutions to non-translation invariant singular SPDEs
 \\ % The assignment title
\horrule{0.5pt} \\[0.5cm] % Thick bottom horizontal rule
}

\author{Harprit Singh\orcidlink{0000-0002-9991-8393}}
\institute{University of Vienna, AT}
\date{\normalsize\today} % Today's date or a custom date

\begin{document}

\maketitle % Print the title
\begin{abstract}
We exhibit a canonical, finite dimensional solution family to certain singular SPDEs of the form 
\begin{equation}
\left(\partial_t-  \sum_{i,j=1}^d a_{i,j}(x,t) \partial_i \partial_j - \sum_{i=1}^d b_i(x,t) \partial_i - c(x,t)\right) u = F(u, \partial u, \xi) \ ,
\end{equation}
where $a_{i,j}, b_i, c: \mathbb{T}^d\times \mathbb{R} \to \mathbb{R}$ 
and
$A=\{a_{i,j}\}_{i,j=1}^d$ is uniformly elliptic.
More specifically, we solve the non-translation invariant g-PAM, $\phi^4_2$, $\phi^4_3$ and KPZ-equation and show that the diverging renormalisation-functions are local functions of $A$. We also establish a continuity result for the solution map with respect to the differential operator for these equations.
%\red{Remove this:}
%\gray{Lastly, we observe that for more singular equations the required renormalisation functions may depend on whether the terms involving $b$, $c$ are interpreted as part of the left or right hand side of the equation.}
%\red{TODO:
%\begin{itemize}
%\item Fix comments along above remark.
%%\item Numerical value cherry for $\phi^4_2$.
%%\item Semi-norm not quite ok, fix as in periodic article.
%%\item More general regularisations: $\rho$ needs volume pre factor.
%\end{itemize}
%}
\end{abstract}
\tableofcontents
%\newpage
\section{Introduction}
While the articles \cite{Hai14}, \cite{BHZ19}, \cite{BCCH20} and \cite{CH16} yield a comprehensive understanding of translation invariant subcritical SPDEs of parabolic type and several results also apply to equations of the form
\begin{equation}\label{eq:type of equations}
  \left(\partial_t-  \sum_{i,j=1}^d a_{i,j}(x,t) \partial_i \partial_j - \sum_{i=1}^d b_i(x,t) \partial_i - c(x,t)\right) u = F(u, \partial u, \xi) \ ,
 \end{equation}
where $a_{i,j}, b_i, c:  \mathbb{T}^d \times\mathbb{R} \to \mathbb{R}$ are smooth functions
and
$A=\{a_{i,j}\}_{i,j=1}^d$ is uniformly elliptic
still, a meaningful solution theory to such singular SPDEs has not been obtained to date. 
%This article aims to fill this gap in the current literature by explicitly renormalising certain important example equations. 
Let us explain where the difficulty lies using the example of the non-linearity $F(u, \partial u, \xi)= -u^3 + \xi$. The solution theory for singular SPDEs involves renormalisation of certain ill-defined products, which in turn involves subtracting (infinite) counter-terms. In the translation invariant setting, this amounts to subtracting diverging constants and since these constants are not unique, one thus does not obtain only one solution, but a finite dimensional family of solutions. This can be rephrased as saying that the natural object of study is not the original (single) equation, but the (finite dimensional) family of equations with right hand sides $F_{C}(u, \partial u, \xi)= -u^3 +Cu  +\xi$ for $C\in \mathbb{R}$. 

When the differential operator is not translation invariant, this redundancy seems to become more severe as one needs to involve renormalisation functions and, at least on first sight, one might end up with an infinite dimensional solution family. 
The main contribution of this article is to present a systematic way to identify a natural \textit{finite dimensional} solution-theory for equations of the form \eqref{eq:type of equations} and to ``solve'' the non-translation invariant analogues of the g-PAM, $\phi^4_2$, $\phi^4_3$ and KPZ-equation. We argue that the solutions obtained in this article are canonical, since the corresponding counter-terms respect the ``geometry'' dictated by the differential operator, see also Remark~\ref{rem:riemannian}, and they are consistent with scaling heuristics when such apply.

\subsection{Scaling heuristics}\label{sec:scaling heuristic}
We first explain how classical scaling heuristics (often) suggest the correct renormalisation functions working with the example 
of the $\phi^4_d$-equation for $d\in \{2,3\}$. Writing $A\cdot \nabla^2:= \sum_{i,j=1}^d a_{i,j} \partial_i \partial_j$ it is formally given by
$$ (\partial_t- A\cdot \nabla^2)u = -u^{3} +\xi \ ,$$
where $\xi$ denotes space time white noise.
Indeed, formally writing $\tilde{u}^\lambda(t,x):=-\lambda^{d/2-1} u(t_0 + \lambda^2 t, x_0 + \lambda x )$ as well as $\tilde{A}^\lambda (t,x)= A(t_0 + \lambda^2 t, x_0 + \lambda x)$ and $\tilde{\xi}^{\lambda}(t,x)=\lambda^{d/2+1} \xi(t_0 + \lambda^2 t, x_0 + \lambda x) $ one finds that
$$(\partial_t-  \tilde{A}^\lambda \cdot \nabla^2) \tilde{u}^\lambda = \lambda^{2-3(d/2-1)} \tilde{u}^3 +\tilde{\xi}^{\lambda} \ .$$
Since $\tilde{\xi}^{\lambda}$ is again a space time-white noise and since $\tilde{A}^\lambda \to A(x_0,t_0)$, as $\lambda\to 0$ this suggests that close to $(x_0, t_0)$
the solution $u$ should be well approximated by the solution
$z^{(x_0, t_0)}$ to 
$$ (\partial_t-  A(x_0, t_0)\cdot \nabla^2)z^{(x_0, t_0)} =  \xi \ .$$
Finally, 
writing $\bar{z}$ for the solution to 
$(\partial_t- \triangle)\bar{z}=  \xi$ and $D(x,t)=|\det(A (x, t))|^{-1/2}$
one formally finds that $\frac{\E[(z^{(x_0, t_0)})^2]}{\E[\bar{z}^2]}= D(x_0, t_0)$
which for $d=2$ suggests to consider the family of equations given by 
\begin{equ}\label{eq:intro_phi^4_2}
(\partial_t- \nabla\cdot A(\cdot ) \nabla)u = -u^3 + 3\alpha^{\<2>}   D(x,t)  u + \xi \ , \qquad \alpha^{\<2>}  \in \mathbb{R}\ .
\end{equ}
For $d=3$ the above argument combined with the main result of \cite{BCCH20}.
%, see in particular Remark~2.23 therein, 
suggests to consider
\begin{equ}\label{eq:intro_phi^4_3}
 (\partial_t- \nabla\cdot A(\cdot ) \nabla)u = -u^3 +\left( 3\alpha^{\<2>}  D(x,t) - 9\alpha^{\<22>} D(x,t)^2 \right) u + \xi \ , \qquad  \alpha^{\<2>} , \alpha^{\<22>}\in \mathbb{R}\ ,
\end{equ}
with the second renormalisation constant $\alpha^{\<22>} $ corresponding to the logarithmic divergence.

The same argument also suggests the correct right hand side for the g-PAM and KPZ equation (where for the former the scaling changes slightly, as one is only working with spatial white noise), see
Theorem~\ref{thm:g-pam} and Theorem~\ref{thm:KPZ}.

%For g-PAM in $d=2$, a similar scaling argument suggests to consider
%$$ (\partial_t- \nabla\cdot A(\cdot ) \nabla)u = f_{ij}(u)\left( \partial_i u \partial_j u - \alpha  C(x,t) a^{i,j}(x,t)  g^2(u)\right) + g(u)\left( \xi  -\alpha' C(x,t) g'(u)\right) \ .$$
%\red{And, finally, for the KPZ-equation}

\begin{remark}
While the simple heuristic just described works for the not ``too singular" equations above, it fails in more singular settings such as for the $\phi^{3}_4$ equation, see Remark~\ref{rem:phi34cherry}, where one should instead Taylor expand the coefficient field.
\end{remark}

\begin{remark}\label{rem:riemannian}
Note that on the torus equipped with a non-standard metric $g=g_{i,j}$ the Laplace-Beltrami operator $\triangle^{g}$ is given by  
$$\triangle^{g} f =g^{i,k} \partial_k\partial_j  f - g^{j,k}\Gamma^l_{jk} \partial_l f \ . $$
Taking that point of view, the fact that in this article one needs to choose renormalisation functions while for the analogue examples in \cite{HS23m} renormalisation constants were sufficient stems from the fact that white noise for the flat metric differs from white noise with respect to the metric $g$ by exactly a multiplicative factor 
$\big(\det (g)\big)^{1/4}$. Thus, the choice of counter-terms advocated here, at least for the g-PAM, $\phi^4_{2}$ and $\phi^4_3$-equation is the unique choice consistent with the intrinsic solution theory in \cite{HS23m}.
Let us here mention the work \cite{Hol10} which considers the renormalisation of quantum field theories on curved spaces, a distinct but related problem. There renormalisation 
ambiguities are reduced by imposing, amongst other axioms, covariance and locality properties on
the operator product expansion.
\end{remark}
\begin{remark}
Let us observe that the exact form of the counterterm as in\footnote{See Theorem~\ref{thm:g-pam}, Theorem~\ref{thm:phi4_3} and Theorem~\ref{thm:KPZ} for precise statements.} \eqref{eq:intro_phi^4_2}\&\eqref{eq:intro_phi^4_3} can only be expected when one works with regularisations which are `sufficiently covariant', such as regularisation using the fundamental solution of the operator in Section~\ref{subsec:heat kernel reg} or the mollifiers of Section~\ref{sec:other regularisation}, see also \cite[Sec.~15]{HS23m}. More general regularisations are discussed in Section~\ref{sec:flat_reg} where it is noted that one can still express the divergent counterterms explicitly and that they are still local functionals of the coefficient field.
\end{remark}

\begin{remark}
Let us also mention the growing body of research studying quasilinear singular SPDEs, 
i.e.~ equations where the differential operator is of the form \newline{$(\partial_t-  \sum_{i,j=1}^d a_{i,j}(u) \partial_i \partial_j )u $.}
By now there exist several approaches to these equations such as \cite{GH19, bru24, bru24a} , \cite{OW19,lin23, ott24,brou25}, 
 \cite{bai24} as well as the works \cite{fur19}, \cite{bai19}.
For these equations one can (and does) ask whether the counterterms can be chosen to be `local' (this time in the solution, or even the coefficient field evaluated at the solution). 
%and whether one can somehow single out a finite dimensional solution family.
\end{remark}

\subsection{Non-translation invariant SPDEs and Regularity Structures}
The analytic machinery of \cite{Hai14} applies to non-translation invariant SPDEs and only second part thereof and the articles \cite{BHZ19} \cite{BCCH20} and \cite{CH16} make significant use of translation invariance. It was shown in \cite{BB21} that preparation maps introduced in \cite{Bru18} are compatible with the use of renormalisation functions.\footnote{The renormalisation functions suggested therein \cite[Eq.~3.15]{BB21} do not correspond to the ones used in this article and, there, convergence of models is not considered.}
The article \cite{HS23m} automates all of the work for a local solution theory to subcritical equations of the form \eqref{eq:type of equations} \textit{except} for the identification of the ``correct'' finite dimensional renormalisation procedure/group, see the discussion in \cite[Section~15]{HS23m}, and the stochastic estimates on models.

\subsubsection*{Structure of the article}
In Section~\ref{sec:heat kernels}, we recall the classical paramerix construction of heat kernels $\Gamma$ associated to the linear part of \eqref{eq:type of equations}, c.f. the bibliographical remarks in \cite{Fri08}. We focus on capturing precise enough asymptotics of $\Gamma$ close to the diagonal, in order to identify  canonical renormalisation functions. In Section~\ref{sec:estimates on counterterms} we check that the constructed kernels fit within the framework of regularity structures and consider the continuity of the abstract solution map with respect to the coefficients of the differential operator.

Then, in Section~\ref{subsec:heat kernel reg} using mollification based on the heat kernel the `correct' renormalisation functions for the g-PAM, $\phi^4$ and the KPZ non-linearity are identified. In particular, we note in Remark~\ref{rem:phi34cherry} that for the more singular $\phi^3_4$ equation the scaling heuristic of Section~\ref{sec:scaling heuristic} does not apply. In Section~\ref{sec:other regularisation} we argue that for the equations treated here more general `covariant' schemes also allow for counter-terms of the same form. Section~\ref{sec:flat_reg} discusses the case of more general regularisations.

Finally, in Section~\ref{sec:main results} we ``solve'' the equations mentioned above. This section is particularly short since the existing literature, specifically \cite{Hai14}, \cite{BB21} and \cite{HS23m}, supplies all necessary results not already established in the previous sections. 

\subsubsection*{Scope for generalisation}
Heat kernels asymptotics are also known for higher order parabolic systems, c.f.\ \cite[Chapter~9]{Fri08} and one expects analogous results to the ones here to hold in that setting, as well as for sub-elliptic/ultra-parabolic equations, c.f.\ \cite[Problem B]{MS25}.
Furthermore, while the arguments in this article allow for weakened regularity assumptions on $A$, c.f.\ Remark~\ref{rem:weakened hölder regularity},
we do not claim that these are optimal.

Lastly, while the methodology presented here robustly applies,
%to parabolic systems of arbitrary degree
it is carried out for each equation by hand. 
This lets one wonder, if it is possible to automate the two steps:
\begin{itemize}
\item Identifying a canonical \textit{finite} dimensional renormalisation group. 
\item Obtaining convergence of the model.
\end{itemize}
These two steps are of somewhat complementary nature, c.f.\ \cite{Hai17, HQ18, CH16, HS23}. 
%and in particular the forthcoming work \cite{Ste23+}.
%The former requires precise knowledge of the heat kernel close to the diagonal, together with arguments reminiscent of negative renormalisation in \cite{HQ18}, \cite{Hai17}, while the latter is less sensitive to upper bounds on heat kernels and derivatives thereof, c.f.~\cite{CH16} and \cite{23}. 
\\
\\
\noindent\textbf{Acknowledgements.}
HS wishes to thank R.~Steele, L.~Broux and M.~Hairer for valuable conversations and M.~Gubinelli additionally for pointing out a misleading formulation in the previous version. HS gratefully acknowledges funding by the Imperial College London President's PhD Scholarship and EPFL Lausanne.

\newpage
\section{Heat Kernels}\label{sec:heat kernels}
In this section we recall the construction of the heat kernel of differential operators on $\mathbb{R}^{d+1}$
of the form
\begin{equation}\label{eq:main_diff_operator}
L = \partial_t - \sum_{i,j} a_{i,j}(x,t) \partial_i \partial_j - \sum_i b_i(x,t) \partial_i - c(x,t)
\end{equation}
following Avner-Friedman \cite[Chapter~1]{Fri08} and collect some of the properties needed here. Without loss of generality, assume that $a_{i,j}=a_{j,i}$.

\begin{assumption}\label{ass:elliptic_and_bounds}
There exists constants $\bar{\lambda}_0, \bar{\lambda}_1>0$ such that for every $z\in \mathbb{R}^{d+1}$ and $\zeta\in \mathbb{R}^d$ we have $$\bar{\lambda}_0|\zeta|^2\lesssim  a_{i,j}(z) \zeta_i \zeta_j \lesssim \bar{\lambda}_1|\zeta|^2 \ $$ 
and the coefficients $a_{i,j}$, $b_i$ and $c$ are smooth and have globally bounded derivatives. 
\end{assumption}

\begin{remark}
Note that Lipschitz continuous coefficients (with respect to the parabolic distance) are sufficient for all considerations of this section, except for the second claim in Proposition~\ref{prop:heat_kernel}, Remark~\ref{rem:upper bound components derivatives}, Lemma~\ref{lem:reflection} and Section~\ref{sec:adjoint} as well as Section~\ref{sec:cont-dependence}.
\end{remark}

\begin{remark}
The argument outlined here is carried out in \cite{Fri08} keeping sharper track of the continuity assumptions on the coefficients. Existence of heat kernels is known under even weaker conditions than treated there, c.f.\ \cite{HK04} and \cite{Aus00}. 
\end{remark}

We shall write $a^{i,j}$ for the entries of the inverse matrix $A^{-1}$ of $A$. We define 
$$\vartheta^{z}(\zeta)= \sum_{i,j} a^{i,j}(z) \zeta_i \zeta _j, \qquad
 w^{z}(x,t; \zeta, \tau)=
 \frac{\mathbf{1}_{ \{ t>\tau \} }}{ (t-\tau)^{d/2}}
  \exp \left(\frac{\vartheta^{z}(x-\zeta)}{4(t-\tau)}\right)$$
and $C(z)= (4\pi)^{-d/2} \det(A(z))^{-1/2}$, where 
$$
\mathbf{1}_{ \{ t>\tau \} }= \begin{cases}
1 & \text{if } t>\tau,\\
0 &\text{else.}
\end{cases}
$$
The fundamental solution of the differential operator with coefficient ``frozen'' at $z=(\zeta, \tau)\in \mathbb{R}^{d+1}$ is given by 
$ C(\zeta,\tau) w^{(\zeta,\tau)}(x,t; y,s)$. We set
\begin{equation}\label{eq:exlicit Z_0}
Z(x,t; \zeta,\tau):= C(\zeta,\tau) w^{(\zeta,\tau)}(x,t; \zeta, \tau)  =  C(\zeta,\tau) w^{(\zeta,\tau)}(x-\zeta,t-\tau; 0, 0)\ .
\end{equation}
One has for any $\lambda^*_0<\lambda_0$, $\mu \in [0, d/2]$ and $T>0$ the estimate
\begin{align}
 |w^{z}(x,t; \zeta, \tau)|=& \left| (t-\tau)^{-\mu} \left( \vartheta^{z}(x-\zeta)\right)^{\mu-d/2}  \left(\frac{\vartheta^{z}(x-\zeta)}{(t-\tau)}\right)^{d/2-\mu}   \exp \left(-\frac{\vartheta^{z}(x-\zeta)}{4(t-\tau)}\right)\right|\nonumber \\
 \lesssim & (t-\tau)^{-\mu} |x-\zeta|^{-d +2\mu}    \exp \left(- \frac{\lambda^*_0 |x-\zeta|^2}{4(t-\tau)} \right)  \ ,\label{eq:upperbound strategy}
\end{align}
uniformly over $(x,t), (\zeta, \tau)\in \mathbb{R}^{d+1}$ satisfying $|t-\zeta|<T$.
One notes the following identity 
\begin{align*}
LZ(x,t; \zeta,\tau)&= \sum_{i,j} (a_{i,j}(x,t)-a_{i,j}(\zeta,\tau)) \partial_i \partial_j Z(x,t; \zeta,\tau) \\
&\qquad +\sum_i \beta_i(x,t) \partial_i Z(x,t; \zeta,\tau) + c(x,t) Z(x,t; \zeta,\tau) 
\ ,
\end{align*}
and thus in particular for any $\mu \in [0, d/2]$, 
\begin{equation}\label{eq:claimed upper}
LZ(x,t; \zeta,\tau)\lesssim (t-\tau)^{-\mu}  |x-\zeta|^{-d-1+2\mu} \exp \left(- \frac{\lambda^*_0 |x-\zeta|^2}{4(t-\tau)} \right) \ .
\end{equation}
Next, for $\nu \geq 1$ recursively define $(LZ)_1:= (LZ)$ and
$$(LZ)_{\nu+1}(x,t; \zeta, \tau) = \int_{\tau}^t \int_{\mathbb{R}^d}   LZ(x,t; \eta, \sigma) (LZ)_\nu (\eta, \sigma; \zeta, \tau) d\eta\, d\sigma \  $$
as well as $Z_0 :=Z$ and
\begin{equation}\label{eq:Z_nu}
Z_\nu (x,t; \zeta, \tau) = \int_{\tau}^t \int_{\mathbb{R}^d} Z(x,t; \eta, \sigma)  (LZ)_\nu (\eta, \sigma; \zeta, \tau) \, d\eta\, d\sigma \ .
\end{equation}
The next lemma collects some properties of these functions and is contained in the discussion around \cite[Eq.~4.14]{Fri08}.
\begin{lemma}\label{lem:bounds_lemma2.4}
For any $\lambda^*_0<\lambda_0$ and $T>0$, there exist positive constants $H_0$, $H$ such that
$$|(LZ)_{\nu}(x,t; \zeta,\tau)|\leq \frac{H_0H^{\nu}}{\Gamma(\nu )} (t-\tau)^{\nu-d/2-1}  \exp\left( -\frac{\lambda^*_0 |x-\zeta|^2}{4(t-\tau)} \right)$$
and 
\begin{equation}\label{eq:upper bound}
|Z_{\nu}(x,t; \zeta,\tau)|\leq \frac{H_0H^{\nu}}{\Gamma(\nu)} (t-\tau)^{\nu-d/2}  \exp\left( -\frac{\lambda^*_0 |x-\zeta|^2}{4(t-\tau)} \right) \ ,
\end{equation}
uniformly over $(x,t), (\zeta, \tau)\in \mathbb{R}^{d+1}$ satisfying $|t-\zeta|<T$.
\end{lemma}

\begin{corollary}\label{cor:bound_on_Z_>nu}
In particular one finds in the setting of the above lemma that for $Z_{\geq \nu}:= \sum_{n \geq \nu}  Z_{ n}$ 
\begin{equation}
|Z_{\geq \nu}(x,t; \zeta,\tau) |\lesssim (t-\tau)^{\nu-d/2}  \exp\left( -\frac{\lambda^*_0 |x-\zeta|^2}{4(t-\tau)} \right)\ .
\end{equation}
\end{corollary}

%%
%%
%%Finally we set 
%%$$\Phi(\eta, \sigma; \zeta, \tau)= \sum_{\nu=1}^\infty (LZ)_\nu (\eta, \sigma; \zeta, \tau) \ .$$
%%\gray{The solution of the integral equation
%%$$
%%\Phi^*(x,t;\zeta, \tau)= L^*Z^*(x,t;\zeta, \tau) + \int_t^\tau \int_{\mathbb{R}^d} L^*Z^*(x,t;\eta, \sigma)\Phi^*(\eta,\sigma;\zeta, \tau)\, d\eta\, d\sigma \ .
%%$$}

The following is the content of \cite[Ch.~1, Thm.~8 \& Ch.~1, Section~6 \& Ch.~9, Sec.~6, Thm.~8]{Fri08} 
\begin{prop}\label{prop:heat_kernel}
Under Assumption~\ref{ass:elliptic_and_bounds} the (unique) fundamental solution of the differential operator \eqref{eq:main_diff_operator}
is given by 
\begin{align*}
\Gamma(x,t; \zeta, \tau)
%&= Z(x,t; \zeta, \tau) + \int_{\tau}^t \int_{\mathbb{R}^d} Z(x,t; \eta, \sigma) \Phi(\eta, \sigma; \zeta, \tau) \, d\eta\, d\sigma\\
&= \sum_{\nu=0}^\infty Z_{\nu}(x,t; \zeta, \tau) \ .
\end{align*}
Furthermore, for $a,a' \in \mathbb{N}^d,\ b,b'\in \mathbb{N}$, $T>0$
$$|\partial_x^{a}\partial_\zeta^{a'} \partial_{t}^b \partial_\tau^{b'} \Gamma(x,t; \zeta, \tau)|\lesssim_{a,a',b, b', T}(t-\tau)^{-(|a|+|a'| + 2b+2b' + d)/2}  \exp\left( -\frac{\lambda^*_0 |x-\zeta|^2}{4(t-\tau)} \right) \ $$
uniformly over $(x,t), (\zeta, \tau)\in \mathbb{R}^{d+1}$ satisfying $|t-\zeta|<T$.
\end{prop}

\begin{remark}\label{rem:upper bound components derivatives}
By proceeding exactly as in the proof of \cite[Ch.~6, Thm.~8]{Fri08} one finds that there exists $H>0$ such that for $a,a' \in \mathbb{N}^d,\ b,b'\in \mathbb{N}$
$$|\partial_x^{a}\partial_\zeta^{a'} \partial_{t}^b \partial_\tau^{b'} Z_{\nu}(x,t; \zeta,\tau)|
\lesssim_{a,a',b, b'}
 \frac{H^{\nu}}{\Gamma(\nu)} (t-\tau)^{\nu-(|a|+|a'| + 2b+2b' + d)/2}   \exp\left( -\frac{\lambda^*_0 |x-\zeta|^2}{4(t-\tau)} \right) \ .$$
 As a consequence the same bound, up to multiplicative constant, holds for $Z_{\geq\nu}$.
\end{remark}

\begin{corollary}
For $\tau<\sigma<t$ the following identity holds
\begin{equation}\label{eq:composition}
\Gamma(x,t;\xi,\tau)= \int_{\mathbb{R}^d} \Gamma(x,t; y,\sigma) \Gamma(y,\sigma;\xi,\tau) \, dy \ .
\end{equation}
\end{corollary}
\begin{proof}
This follows directly from the uniqueness of the fundamental solution.
\end{proof}
%The next lemma follows by a simple computation.
%\begin{lemma}\label{lem:derivative_estimate}
%One has the following bounds
%\begin{align*}
%|\partial_{i;x} \Gamma(\zeta+ x,t;\zeta,\tau)|&\lesssim \frac{1}{(t-\tau)^{d/2}}   \exp \left(- \frac{\lambda^*_0 |x-\zeta|^2}{4(t-\tau)} \right)\\
%|\partial_{t} \Gamma(x,t+\tau;\zeta,\tau)|&\lesssim \frac{1}{(t-\tau)^{d/2}}   \exp \left(- \frac{\lambda^*_0 |x-\zeta|^2}{4(t-\tau)} \right)
%\end{align*}
%\end{lemma}\harprit{double check/ cite}
%
\subsection{Behaviour close to the diagonal}
We write 
\begin{equation}\label{eq:reflection}
\mathfrak{R}_\zeta: \mathbb{R}^d\to \mathbb{R}^d,\qquad  x\mapsto -x+2\zeta
\end{equation}
 for the map reflecting at the point $\zeta\in \mathbb{R}^d$ and observe that 
$$Z_0(\mathfrak{R}_\zeta(x),t; \zeta, \tau)=Z_0(x,t; \zeta, \tau)\ . $$
 
 \begin{lemma}\label{lem:reflection}
One can write $Z_1(z; \bar{z})= \bar{Z}_1(z; \bar{z}) +R_1(z; \bar{z})$ such that $\bar{Z}_1$ satisfies
 \begin{equation}\label{eq:reflection_symmetry}
\bar{Z_1}(\mathfrak{R}_\zeta(x),t; \zeta, \tau)=-\bar{Z}_1(x,t; \zeta, \tau)  ,
 \end{equation}
 as well as the same upper bound \eqref{eq:upper bound} as $Z_1$ (with possibly a larger constant)
and $R_1$ satisfies the upper bound \eqref{eq:upper bound} for $\nu=2$ (with possibly a larger constant).
\end{lemma}
 
% 
% \begin{equation}\label{eq:reflection_symmetry}
%Z_1^{i,j}(\mathfrak{R}_\zeta(x),t; \zeta, \tau)=-Z_1^{i,j}(x,t; \zeta, \tau)  \qquad  Z_1^{i}(\mathfrak{R}_\zeta(x),t; \zeta, \tau)= -Z_1^{i}(x,t; \zeta, \tau)\ .
% \end{equation}
% 
%% 
%
%\begin{lemma}
%Next we we want to decompose $Z_1(z; \bar{z})= \bar{Z}_1(z; \bar{z}) +R_1(z; \bar{z})$ such that $\bar{Z}_1$ satisfies the analogue of \eqref{eq:scaling_lowest_order} with $\lambda^{-|\fraks|+2}$ replaced with $\lambda^{-|\fraks|+3}$ 
%\end{lemma}

%and 
%$R(z; \bar{z})= \mathcal{O}(|z-\bar{z}|)$ for $z\to \bar{z}$. 
\begin{proof}
We find that 
\begin{align*}
Z_1(x,t; \zeta, \tau)=&\int_{\tau}^t \int_{\mathbb{R}^d} Z(x,t; \eta, \sigma)  (LZ)_1 (\eta, \sigma; \zeta, \tau) \, d\eta\, d\sigma\\
=&\sum_{i,j} 
\int_{\tau}^t \int_{\mathbb{R}^d} Z(x,t; \eta, \sigma)
 (a_{i,j}(\eta, \sigma)-a_{i,j}(\zeta,\tau)) \partial_i \partial_j Z(\eta, \sigma; \zeta,\tau) \, d\eta\, d\sigma\\
 &
+\sum_i  
\int_{\tau}^t \int_{\mathbb{R}^d} Z(x,t; \eta, \sigma) \beta_i(\eta, \sigma) \partial_i Z(\eta, \sigma; \zeta,\tau)\, d\eta\, d\sigma\\
& + \int_{\tau}^t \int_{\mathbb{R}^d} Z(x,t; \eta, \sigma)c(\eta, \sigma) Z(\eta, \sigma; \zeta,\tau) \, d\eta\, d\sigma\\
=&\bar{Z}_1(x,t; \zeta, \tau)+ R_1 (x,t; \zeta, \tau) 
\end{align*}
where $\bar{Z}_1:= \sum_{i,j} \bar{Z}_1^{i,j}+\sum_i  \bar{Z}_1^{i}$ for 
\begin{align*}
&\bar{Z}_1^{i,j}(x,t; \zeta, \tau):= \sum_k  \partial_k a_{i,j}(\zeta,\tau) C(\zeta, \tau) \int_{\tau}^t \int_{\mathbb{R}^d} w^{\zeta,\tau}(x,t; \eta, \sigma)
  (\eta-\zeta)_k \partial_i \partial_j Z(\eta, \sigma; \zeta,\tau) \, d\eta\, d\sigma\\
&  \bar{Z}_1^{i}(x,t; \zeta, \tau):=   b_{i}(\zeta,\tau) C(\zeta, \tau) \int_{\tau}^t \int_{\mathbb{R}^d} w^{\zeta,\tau}(x,t; \eta, \sigma)
  \partial_i  Z(\eta, \sigma; \zeta,\tau) \, d\eta\, d\sigma\ ,
\end{align*}
and 
 $$R_1= \sum_{i,j} R^{i,j} +\sum_{i} R^i + Z_1^{0}$$
for
 $R^{i,j}= (Z^{i,j}-\bar{Z}^{i,j})$, $R^{i}= (Z^{i}-\bar{Z}^{i})$ and $$Z_1^{0}(x,t; \zeta, \tau):= \int_{\tau}^t \int_{\mathbb{R}^d} Z(x,t; \eta, \sigma)c(\eta, \sigma) Z(\eta, \sigma; \zeta,\tau) \, d\eta\, d\sigma \ .$$
At this point one sees that $\bar{Z}_1$ and $R_1$ satisfy the properties claimed in the lemma.
\end{proof} 

It is sometimes useful to work with modifications of the expansion in Prop.~\ref{prop:heat_kernel}, for example in view of the above lemma with
\begin{align}
\Gamma(x,t; \zeta, \tau)
%&= Z(x,t; \zeta, \tau) + \int_{\tau}^t \int_{\mathbb{R}^d} Z(x,t; \eta, \sigma) \Phi(\eta, \sigma; \zeta, \tau) \, d\eta\, d\sigma\\
&=  \sum_{\nu\geq 0}^\infty \tilde{Z}_{\nu}(x,t; \zeta, \tau) \ .
\end{align}
where  $\tilde{Z}_{1}= \bar{Z}_{1}$ and 
$\tilde{Z}_{2}= {Z}_{2} + R_1$ while
 $\tilde{Z}_{\nu}(x,t; \zeta, \tau)= {Z}_{\nu}$ for $\nu \notin \{1,2\}$.

\subsection{Adjoint Kernels}\label{sec:adjoint}

The (formal) adjoint operator of $L$ is given by 
$$L^*= \partial_t + \sum_{i,j} a_{i,j}(\cdot) \partial_i \partial_j + \sum_i b^*(x,t)_i \partial_i + c^*(x,t)\ , $$
where $b^*_i =-b_i+ 2 \sum_j \partial_j a_{i,j}$ and $c^*=c -  \sum_i \partial_i b_{i}+\sum_{i,j} \partial_i\partial_j a_{i,j}$.
Its fundamental solution (back-wards in time, i.e. for $t<\tau$) can be constructed in essentially the same way as the one for $L$ and is then given by
$$\Gamma^*(x,t;\zeta, \tau)= \sum_{\nu=0}^\infty Z_{\nu}^*(x,t;\zeta, \tau) $$
where 
$$Z_0^*(x,t;\zeta, \tau)= 
\frac{C(\zeta, \tau) \mathbf{1}_{\{\tau> t  \}}}{(\tau-t)^{d/2}}
  \exp\left( \frac{\vartheta^{(\zeta,\tau)}(x-\zeta) }{4(\tau-t)} \right) $$
and $Z_{\nu}^*(x,t;\zeta, \tau)$ is defined analogously to \eqref{eq:Z_nu} for $\nu\geq 1$. It is the content of \cite[Thm.~15]{Fri08} that one has the identity 
\begin{equation}\label{eq:adjoint}
\Gamma(\zeta, \tau; x,t)=\Gamma^*(x,t;\zeta, \tau) \ .
\end{equation}
For later use we define for $i\in \{1,...,d\}$ 
$$Z_{0;i} ( \eta, \tau; x,t ) =
 -\frac{\sum_l a^{i,l}(\eta-x)_l}{2(t-\tau )} Z_0^*( \eta, \tau; x,t ) \ .$$
\begin{lemma}\label{lem:bound on Z_{0;1}}
One can write $\partial_{x_i} Z_0^*( \eta, \tau; x,t )= Z_{0;i} ( \eta, \tau; x,t ) +R^*_{0,i}( \eta, \tau; x,t )$ where
$$|R^*_{0,i}( \eta, \tau; x,t )|\lesssim_T (t-\tau)^{-d/2}  \exp\left( -\frac{\lambda^*_0 |x-\eta|^2}{4(t-\tau)} \right) \ ,$$
uniformly over $(x,t), (\eta, \tau)\in \mathbb{R}^{d+1}$ satisfying $|t-\zeta|<T$.
\end{lemma}
\begin{proof}
One finds that
$$\partial_{x_i} Z_0^*( \eta, \tau; x,t )= -\frac{\sum_l (a^{i,l}+ a^{l,i})(\eta-x)_l +\langle(\eta-x), (\partial_{x_i}A(x,t))(\eta-x)\rangle }{4(t-\tau)} Z_0^*( \eta, \tau; x,t ), $$
thus the claim follows by observing that
$$Z_{0;i} ( \eta, \tau; x,t ) = -\frac{\sum_l (a^{i,l}+ a^{l,i})(\eta-x)_l}{4(t-\tau)} Z_0^*( \eta, \tau; x,t )$$
together with an estimate as in \eqref{eq:upperbound strategy} .
\end{proof}

%
% \red{we see that
%$$\partial_i Z_0^*( \eta, \tau; x,t )= -\frac{\sum_l (a^{i,l}+ a^{l,i})(\eta-x)_l +\langle(\eta-x), (D_xA(x,t))(\eta-x)\rangle }{4(t-\tau)} Z_0^*( \eta, \tau; x,t ) $$
%and that by Remark~\ref{rem:upper bound components derivatives} any contribution of $Z_{\geq 1}^*$ to $E\left[\bfPi^{\epsilon} \<b2>_{ij} (x,t)\right] $ converges to a finite limit.
%
%We further write $\partial_i Z_0^*= Z_{0;i}^* +R^*_{0,i}$ where
%\begin{equation}\label{eq:decomposition partial_i Z_0^*}
%Z_{0;i} ( \eta, \tau; x,t ) = -\frac{\sum_l (a^{i,l}+ a^{l,i})(\eta-x)_l}{4(t-\tau)} Z_0^*( \eta, \tau; x,t )=
% -\frac{\sum_l a^{i,l}(\eta-x)_l}{2(t-\tau )} Z_0^*( \eta, \tau; x,t ) \ .
% \end{equation}
%and $|R^*_{0,i}|$ is directly seen to be bounded by the right hand side of \eqref{eq:upper bound} with $\nu=0$.}

\subsection{Operators depending on a parameter}\label{sec:cont-dependence}
Consider the situation where $a^{\lambda}_{i,j}(x,t)$, $b^{\lambda}_i(x,t)$ and $c^{\lambda}(x,t)$ depend on a parameter $\lambda\in [0,1]$ and consider the fundamental solutions
$\Gamma^\lambda= \sum_{\nu \geq 0} Z_\nu^\lambda$ as in Proposition~\ref{prop:heat_kernel} to the differential operator
$$L^{\lambda} = \partial_t - \sum_{i,j} a^{\lambda}_{i,j}(x,t) \partial_i \partial_j - \sum_i b^{\lambda}_i(x,t) \partial_i - c^{\lambda}(x,t) \ .$$
The following proposition is well known and follows as direct consequence of the construction in the previous section, c.f.\ \cite[Theorem~2.48]{BGV03} for a similar result in a slightly different setting.
\begin{prop}\label{prop continuity}
Let $A^\lambda=\{a^{\lambda}_{i,j}\}$ satisfy Assumption~\ref{ass:elliptic_and_bounds} for each $\lambda\in [0,1]$. Assuming furthermore that for each $(x_0,t_0)\in \mathbb{R}^{d+1}$ the map
 $\lambda \mapsto (a^{\lambda}_{i,j}(x_0,t_0),  b^{\lambda}_i(x_0,t_0), c^{\lambda}(x_0,t_0))$ is continuously differentiable and that the functions
 $\frac{\partial}{\partial\lambda}  a^{\lambda}_{i,j}, \  \frac{\partial}{\partial\lambda} b^{\lambda}_i,\ \frac{\partial}{\partial\lambda} c^{\lambda}$ are bounded with bounded derivatives (with respect to $(x,t)\in \mathbb{R}^{d+1}$), it holds that for each $z\neq \bar{z}$ the maps
$$ \lambda \mapsto \Gamma^\lambda (z,\bar{z}), \qquad \lambda \mapsto Z^\lambda_{\nu} (z,\bar{z})$$ are continuously differentiable for each $\nu\geq 1$ and the kernels
$\frac{\partial}{\partial\lambda}  \Gamma^\lambda, \ \frac{\partial}{\partial\lambda}Z^\lambda_{\nu} $ and their derivatives (in $z, \bar{z}$)
satisfy identical upper bounds (up to multiplicative constant) as 
$\Gamma^\lambda$ respectively $Z^\lambda_{\nu}$.
\end{prop}

\begin{proof}
One first checks differentiability in $\lambda$ and upper bounds on 
$Z^\lambda$ and $L^\lambda Z^\lambda$ (including \eqref{eq:upperbound strategy} and \eqref{eq:claimed upper}). Then, the stated claims on $Z^\lambda_\nu$ follow by direct computation as for Lemma~\ref{lem:bounds_lemma2.4}. Thus, by summability the claims about $\Gamma^\lambda$ follow as well.
\end{proof}

\subsection{Kernels in Regularity Structures}\label{sec:estimates on counterterms}

The aim of this section is to check that the kernels constructed above fit into the framework of regularity structures and to establish continuity bounds with respect to the kernels.
Thus, throughout this section we assume familiarity with the content of the of the article \cite{Hai14} up to and including Section~7.

First, we introduce some notation in order to quantify
the qualitative \cite[Assumption~5.1]{Hai14}. For $\gamma>0$ and $f: \mathbb{R}^{d}\to \mathbb{R}$ compactly supported such that $D^{k}f$ exists whenever $|k|_\fraks<\gamma$, let
$$ \opP_{z_{0}}^\gamma [f](z):= \sum_{|k|_\fraks<\gamma} \frac{D^{k}f(z_0)}{k!} (z-z_0)^k . $$ 
Similarly, for $\gamma'>0$ and a compactly supported function of two variable $F: \mathbb{R}^{d}\times \mathbb{R}^d \to \mathbb{R}$ such that
$D^{k_1}D^{k_2}F$ exists whenever $|k_1|_\fraks<\gamma$, $|k_2|_\fraks<\gamma'$,
 let
$$\opP^{(\gamma,\gamma')}_{({z}_0, \bar{z}_0)}[F](z, \bar{z}) := \sum_{|k|<\gamma, |l|<\gamma' } \frac{D_1^{k} D_2^{l} F({z}_0, \bar{z}_0) }{k!\ l!} (z-z_0)^k(\bar{z}-\bar{z}_0)^l \ .$$

\begin{definition}
Let $K: \mathbb{R}^{d+1}\times \mathbb{R}^{d+1}\setminus \triangle \to \mathbb{R}$ be a kernel
which can be decomposed as $K(z,z')= \sum_{n\geq 0} K_n(z,z')$ where each $K_n$ is supported on 
$ \{(z,z')\ : \ \|z-z'\|_\fraks \leq 2^{-n+1} \} .$ 
On the space of such kernels $\pmb{\mathcal{K}}$, we define for $L, R \in \mathbb{R}_+$, 
the norm $\|K\|_{\beta;L,R}$ as the smallest number $C$ such that there exists a decomposition $K= \sum K_n$ for which the following bounds are satisfied:
\begin{itemize}
\item For any multi-indices  $k_1, k_2$ satisfying $|k_1|_\fraks < L$, $|k_2|_\fraks< R$ and $n \in \mathbb{N}$ it holds that
\begin{equation}\label{upper bounds}
 \sup_{z,z'}|D_1^{k_1}D_2^{k_2} K_n(z,z') | \leq C2^{(|\fraks|-\beta +|k_1|_\fraks +|k_2|_\fraks )n}
\end{equation}
%uniformly over $n\geq 0$ and $z,z'\in \mathbb{R}^{d+1}$.
as well as
\begin{align*}
& \sup_{z,z',\bar{z}'} \frac{|D^{k_1}_1K(z, z')- \opP^{R}_{\bar{z}' }[D^{k_1}_1K(z,\cdot ) ](\bar{z}')|  }{|z'-\bar{z}'|^R } \lesssim C2^{n(|\fraks|-\beta +|k_1|_\fraks+ R  )}\\
&  \sup_{z,z',\bar{z}} \frac{|D^{k_2}_1K(z, z')- \opP^{R}_{\bar{z}' }[D^{k_2}_2K(\cdot,z' ) ](\bar{z})|  }{|z-\bar{z}|^L } \lesssim C2^{n(|\fraks|-\beta +|k_2|_\fraks+ L  )}
\end{align*}
and 
\begin{equ}
 \sup_{z,\bar{z},z',\bar{z}'}  \frac{|K_n(z,z') -\opP^L_{\bar{z}} [K_n(\,\cdot\, , {z}')](z) -\opP^{R}_{\bar{z}'}[K_n({z}, \,\cdot\,  )](z')+  \opP^{(L,R)}_{(\bar{z}, \bar{z}')}[K_n](z, {z}')| }{|z-\bar{z}|_\fraks^\gamma  |z'-\bar{z}'|_\fraks^{\gamma'} } 
 \leq C2^{n(|\fraks|-\beta +L+ R  )}
\end{equ}
uniformly over $n\geq 0$.
\item Let $\mathbf{k}_{<R}:= \big\{k \in \mathbb{N}^d \ : \ |k|_{\fraks}<R \big\}$ and denote by $\partial\mathbf{k}_{<R}:= \{
k\notin \mathbf{k}_{<R}\ : k-e_{\min\{ j\ : k_j\neq 0  \}\in \mathbf{k}_{<R} }
\}
$ its boundary in the sense of \cite[App.~A]{Hai14}. Then, for any $k_1\in  \mathbf{k}_{<R}$ and $k_2\in \mathbf{k}_{<R}\cup\partial\mathbf{k}_{<R}$
\begin{equation}\label{eq:technical bound}
|\int_{\mathbb{R}^d} (z-z')^{k_1} D^{k_2}_2 K_n(z,z') dz| \leq C2^{-\beta n }
\end{equation}
uniformly over $n\geq 0$ and $z'\in \mathbb{R}^{d+1}$.
\end{itemize}
Finally, we define the vector space $\bfK^\beta_{L,R}= \{ K\in \bfK \ : \|K\|_{\beta;L,R}<+\infty\}$ equipped with the norm $\|\cdot\|_{\beta;L,R}$.
\end{definition}

We observe that the first item of the above definition quantifies \cite[Eq.~(5.4)]{Hai14} and the second item \cite[Eq.~(5.5)]{Hai14} for which we also used that only the terms where derivatives fall on the second variable are actually need (an observation already made in \cite[Sec.~7]{HS23m}).
Note furthermore that, while in \cite[Assumption~5.4]{Hai14} it is assumed that the kernels $K_n$ annihilate polynomials up to a certain degree, this assumption can be replaced by the following assumption which is more convenient in the non-translation invariant setting and can always be satisfied by replacing $\beta$ by $\beta-\epsilon$ for (certain) arbitrarily small $\epsilon\in (0,1)$.
\begin{assumption}\label{ass:irrational_beta}
For all $\alpha\in A$ (the index set of the structure space $T= \bigoplus_{\alpha\in A} T_\alpha$), one has $\alpha+\beta \notin \mathbb{N}$.
\end{assumption}

In what follows we observe that the general abstract fixed point theorem of \cite{Hai14} is continuous with respect to the norms $\|\cdot \|_{\beta;L,R}$ provided $L,R$ are large enough. Since this only amounts to a careful reading of \cite{Hai14} and a more sophisticated variant of these ideas was carried out in \cite{GH19}, we shall only sketch the argument.

Recall, that in the theory of regularity structures one lifts a $K$ to an abstract operator $\mathcal{K}$. This involves, for a given sector $V\subset T$, three operators 
$${I}: V\to T ,\qquad {J}^K(z): V \to \bar{T}, \qquad \mathcal{N}^K_\gamma: \mathcal{D}^\gamma (V) \to C(\mathbb{R}^{d+1},\bar T)\ ,$$
see \cite[Eq.~5.11, Eq.~15 \& Eq.~5.16]{Hai14}.
One readily checks that if the sector $V$ satisfies for $\alpha, \eta \in \mathbb{R}$ that 
$$V\subset \bigoplus_{\alpha\leq \delta \leq \eta} T_\delta\ , $$ the following holds:
\begin{itemize}
\item The linear map $\bfK^\beta_{L,R} \ni K\mapsto J^K\in C( \mathbb{R}^{d+1}, L(V,{T}))$ is Lipschitz continuous if
\begin{equation}\label{LR-assumption}
-R< \alpha, \qquad L>\eta +\beta\ .
\end{equation}
\item The map $\bfK_{L,R}^\beta \ni K\mapsto\mathcal{N}^K_\gamma\in L(\mathcal{D}^\gamma(V),C(\mathbb{R}^{d+1},\bar T))$ is Lipschitz continuous if
\begin{equation}\label{LR-assumption2}
-R< \alpha, \qquad L>\gamma +\beta\ .
\end{equation}
\end{itemize}
Importantly, by a careful reading of the corresponding sections of \cite{Hai14} one finds that for $L,R>0$ satisfying \eqref{LR-assumption}\& \eqref{LR-assumption2}, the results \cite[Thm~5.12, Thm.~5.14, Prop.~6.16 and Theorem~7.1]{Hai14} still hold in the slightly more general setting where $Z=(\Pi, \Gamma)$ realises a kernel $K\in \bfK^\beta_{L,R}$ for $I$ and $\bar Z=(\bar \Pi, \bar \Gamma)$ realises a distinct kernel $\bar K\in \bfK^\beta_{L,R}$ for $I$, but with the implicit constant in the inequalities of these results depending on a constant  $C>0$ for which the inequality then holds
uniformly in kernels $K,\bar{K}$ satisfying $\|K\|_{\beta;L,R}\vee \|\bar K \|_{\beta;L,R}<C$ and 
with an additional term $\|K-\bar{K} \|_{\beta;L,R}$ on the respective right hand sides.

For a regularity structure equipped with an abstract integration map $I: V\to T$ we define 
$$\bfK_{L,R}^\beta\ltimes\mathcal{M}\subset \bfK_{L,R}^\beta\times\mathcal{M}$$ to consist of pairs $(K,Z)$ such that 
 $Z$ realises $K$ for $I$. 
%\gray{
%Recall that the abstract lift of the fix point problem in \cite{Hai14} also involves the abstract lift of a smooth compactly supported kernel $K_{-1}$. 
%Here, we define the space $C^{L,R}$ to consist of continuous functions $K_{-1}: \mathbb{R}^{d+1}\times \mathbb{R}^{d+1} \to \mathbb{R}$ which are supported within finite distance from the diagonal and such that the following norm is finite
%\begin{align*}
%\|K_{-1}\|_{C^{L,R}}=& \max_{|k_1|<L, |k_2|<R} \sup_{z,z'\in \mathbb{R}^{d+1}} |D_1^{k_1}D_2^{k_2} K_n(z,z') |\\
%& +  \max_{k\in \mathbb{M}_L, k_2\in \mathbb{M}_R}\sup_{z,z'\in \mathbb{R}^{d+1}}\sup_{|h_1|, |h_2|<1} \frac{|D_1^{k_1}D_2^{k_2} K_n(z+h_1,z'+h_2)-D_1^{k_1}D_2^{k_2} K_n(z,z')|}{|h_1|_\fraks^{L-|k_1|_\fraks}+ |h_2|_\fraks^{R-|k_2|_\fraks} } \ ,
%\end{align*}
%where for $C>0$ we donote by $\mathbb{M}_C\subset \mathbb{N}^{d+1}$ the maximal elements of $\{k\in \mathbb{N}^{d+1} \ : \ |k|_\fraks<C\}$.
%We shall use the notation 
%$$\opP^\gamma: \mcC^\gamma (\mathbb{R}^{d+1}) \to \mathcal{D}^\gamma (\bar{T})\, \qquad f\mapsto \opP^\gamma[f]= \sum_{|k|_\fraks<\gamma} \frac{X^{k}}{k!}  D^{k}f,$$
%then we can lift a kernel $K_{-1}\in C^{L,R}$  to modelled distributions exactly as in \cite[Eq.~7.7]{Hai14} by
%$$\mathcal{D}^\delta (V) \mapsto  \mathcal{D}^\gamma (\bar{T}), \qquad f\mapsto \opP^\gamma[ {K}_{-1}\mathcal{R} f]\ ,$$
%where $\mathcal{R}$ denotes the reconstruction operator.
%}
One concludes that the following \textit{slight} generalisation of \cite[Theorem~7.8]{Hai14} holds.
\begin{theorem}\label{thm:fixed point}
In the setting of \cite[Theorem~7.8]{Hai14}, under the Assumption~\ref{ass:irrational_beta} instead of \cite[Assumption~5.4]{Hai14} and assuming 
that $\gamma<L$ and that $-R$ is smaller than the regularity of the sector $V$, consider 
the solution map for $T>0$ small enough to the fixed point problem
$$u= (I + J^{K} + \mathcal{N}^K) \mathbf{R}^+ F(u)  + v\ ,$$
 as a map $$S_T:  \left (\bfK_{L,R}^\beta\ltimes\mathcal{M}\right) \times \mathcal{D}^{\gamma,\eta}_P(\bar T) \to \mathcal{D}^{\gamma,\eta}_P
 , \qquad ((K,Z), v)\mapsto u \ .$$
 
Then, assuming $F$ is strongly locally Lipschitz, for each bounded\footnote{
with respect to the distance 
$$
 \left(((K,Z), v), ((\bar K, \bar Z), \bar v)\right)  
 \mapsto 
 \|K-\bar{K}\|_{\beta;L,R} 
  + \vertiii{Z,\bar{Z}}_{\gamma;O} 
%   +
%  \|K_{-1}-\bar{K}_{-1}\|_{C^{L,R}} 
  +\vertiii{v;\bar{v}}_{\gamma,\eta; T} \ , $$} 
  subset \linebreak $B\subset \left (\bfK_{L,R}^\beta\ltimes\mathcal{M}\right)\times \mathcal{D}^{\gamma,\eta}_P(\bar T)$
there exists $T>0$ such that the solution map $S_T$ is well defined and Lipschitz on $B$.
%$$
% \left(((K,Z), K_{-1},v), ((\bar K, \bar Z), \bar{K}_{-1},\bar v)\right)  
% \mapsto 
% \|K-\bar{K}\|_{L,R} 
%  + \vertiii{Z,\bar{Z}}_{\gamma;O} 
%   +
%  \|K_{-1}-\bar{K}_{-1}\|_{C^{L,R}} 
%  +\vertiii{v;\bar{v}}_{\gamma,\eta; T} \ , $$
%on  $\left (\bfK_{L,R}^\beta\ltimes\mathcal{M}\right)\times C^{L,R} \times \mathcal{D}^{\gamma,\eta}_P(\bar T)$.
\end{theorem}

\begin{remark}\label{rem: kernel LR assumption for specific equations}
Let us consider what choices of $L,R$ in the above Theorem~\ref{thm:fixed point} work for our example equations. 
\begin{enumerate}
\item For g-PAM in $d=2$ one can take $R>1, L>1$, c.f.\ \cite[Sec.~9.1\& Sec.~9.3]{Hai14}.
\item For the $\phi^4_2$-equation one can take $R>2, L>0$. But one can observe that after making sense of $\<1>$ ``by hand'' instead of using the extension theorem of \cite[Thm.~5.14]{Hai14}, it suffices to take $R,L>0$, c.f. \cite{DD03, TW18, HSper}. 
\item For the $\phi^4_3$-equation, after making sense of $\<1>$ and $\<30>$ ``by hand'', one can take $R>1$ and $L>1$, c.f.\ \cite[Sec.~9.2 \& Sec.~9.4]{Hai14} and \cite{CMW23}.
\item Following the usual convention for the KPZ-equation of writing $\<I>$ for the derivative of the heat kernel and $\<bI>$ for the heat kernel, after making sense of $\<b1>$, $\<20b>$ and $\<210b>$ ``by hand'' one can take $R>1/2$ and $L>3/2$ for this equation, c.f.\ \cite{Hai13} and \cite{FH20}.
%\gray{
%$$\mathfrak{T}^\mathcal{J}:=\left\{ \<1>, \; \<20>, \;\<210>
%\;\<1'>, \; \<2'0>\ \right\}$$
%and
%$${\mathfrak{T}}^r:= 
%\left\{\ \<X>, \; \<2>, \; \<21>,  \; \<211>, \; \<22j>,
%\; \<2'>, \; \<2''>, \; \<2'1> ,\; \<21'>\ \right\}\ . $$
%$L=2$ and $R=2$, after making sense of $\<1>$ and  $\<20>$ by hand maybe $R=1$.}
\end{enumerate}
\end{remark}
\begin{remark}
Note that we could have weakened the condition 
\eqref{eq:technical bound} slightly further by imposing only a H\"older type estimate instead of allowing $k_2\in \partial\mathbf{k}_{<R}$. 
\end{remark}

Finally, we check that non-translation invariant kernels constructed here fit within the framework of regularity structures. From now on, for the rest of this article we make the following assumption: 
\begin{assumption}\label{ass:periodic}
The functions $a_{i,j}(\ \cdot \ ,t), b_i(\ \cdot \ ,t), c(\ \cdot \ ,t)$ are periodic for every $t\in \mathbb{R}$.
\end{assumption}

 We work with a truncation of the heat kernel in time specified by $\kappa: \mathbb{R}\to [0,1]$ such that 
 \begin{itemize}
 \item $\kappa(t)=0$ for $t<0$ or $t>1$,
 \item $\kappa(t)=1$ for $t\in (0,1/2)$.
 \item $\kappa|_{\mathbb{R}_+}$ is smooth.
 \end{itemize}
 
We also fix $\phi: \mathbb{R}^d\to [0,1]$ smooth and compactly supported on ${B}_{1+1/10}(0)\subset \mathbb{R}^d$ such that $\sum_{k\in \mathbb{Z}^d} \phi(x+k)=1$ for all $x$.

\begin{prop}\label{lem:fits in reg struct}
Under Assumption~\ref{ass:elliptic_and_bounds} and Assumption~\ref{ass:periodic} the kernel 
$$K(x,t, \zeta,s)=\sum_{k\in \mathbb{Z}^d}\kappa(t-s)\phi(x-\zeta)\Gamma(x,t,\zeta+k,s)$$ is non-anticipative and satisfies 
$\|K\|_{\beta; L,R}< +\infty$
for any $L,R\geq 0$ and $\beta\leq 2$. 
Furthermore, in the setting of Proposition~\ref{prop continuity} the family of kernels
$$K^\lambda=\sum_{k\in \mathbb{Z}^d}\kappa(t-s)\phi(x-\zeta)\Gamma^\lambda(x,t,\zeta+k,s)$$
additionally satisfies
$
\|K^\lambda-K^\mu \|_{\beta; L,R}\lesssim |\lambda-\mu|
$
 for every $L,R\geq 0$ and $\beta\leq 2$. 
\end{prop}

\begin{proof}
It is clear that $K$ is non-anticipative.
Let $\varphi\in \mathcal{C}^{\infty}_{c}(B_2\setminus B_{1/2})$ be such that \newline$\sum_{n=0}^{\infty}\varphi\big(\mathcal{S}_\fraks^{2^{-n}}(z) \big)=1$ for every $z\in B_{1+1/10}\setminus\{0\}$ and  set 
$
{K}_n(z;z')=  \varphi \big(\mathcal{S}_\fraks^{2^{-n}}(z-z') \big)		K(z;z')  \ 
$ for $n\geq 0$.
In order to see that $\|K\|_{\beta; L,R}< +\infty$, only the bound \ref{eq:technical bound} requires an argument, the rest follow straightforwardly from Proposition~\ref{prop:heat_kernel}. We first write $$\Gamma= \sum_{\nu=0}^{\lfloor R \rfloor} Z_\nu + Z_{>\lfloor R \rfloor}\ .$$
For the summand
$Z_{>\lfloor R \rfloor}$ the bound follows directly form Remark~\ref{rem:upper bound components derivatives}. For the the summand $Z_0$, the bound follows by first turning the derivative $D_2$ into terms only involving derivatives $D_1$ in the first variable using the explicit formula~\eqref{eq:exlicit Z_0} and then integrating by parts. 
The bound on 
$Z_\nu$ for $\nu \in \{1,...,\lfloor R \rfloor\}$ follows similarly, but by additionally using \eqref{eq:Z_nu} together with the explicit form of $LZ$ above \eqref{eq:claimed upper}.
The second part of the lemma follows analogously, but furthermore using Proposition~\ref{prop continuity}.
\end{proof}

\begin{remark}
Note that in \cite[Sec.~7]{Hai14} the kernel is decomposed into a compactly supported singular part and a compactly supported smooth part,
where the singular part is supported on a fundamental domain. This last property is convenient in establishing convergence of models, but not crucial.
\end{remark}

\begin{remark}\label{rem:weakened hölder regularity}
Note that by keeping slightly more careful track of the regularity assumptions on $A$, one finds that Assumption~\ref{ass:elliptic_and_bounds} can be weakened depending on the equation at hand. For example, one finds in view of Remark~\ref{rem: kernel LR assumption for specific equations}
that for the $\Phi^4_2$ equation one can weaken the regularity assumption on $A$ to mere H\"older continuity for some positive exponent using the bounds in \cite[Chapter~1, Section~6]{Fri08}.
\end{remark}

\section{Canonical Counterterms}

Finally, we turn to the main contribution of this article, the identification of canonical counterterms/solutions to the non-translation invariant g-PAM, $\phi^4_2$, $\phi^4_3$ and KPZ equation. We shall freely use standard graphical notations such as $\<1>= I\Xi$ and $\<2>,\<22>$, c.f.\ \cite{FH20}, \cite{HQ18} and \cite{HP15}; the meaning will be clear from the specified context.

%\gray{
%
%
%We shall always fix $\rho\in C^\infty(\mathbb{R}^d+1)$ a radially symmetric smooth test function and write $\rho_{\epsilon}$ for the rescaled test function $\rho_{\epsilon}(z)= \frac{1}{\epsilon^{|\fraks|}}\rho_{\epsilon}(\mathcal{S}_\epsilon z)$.}

%We shall freely use standard graphical notations such as $\<1>= I\Xi$ and $\<2>,\<22>$ where the meaning will be clear from the context of the equation considered.
%In particular, here $I$ will be the abstract lift of $K$ respectively $\nabla K$, where we only work with truncations of the kernel in time.

Throughout this section, $\alpha$ will denote a constant, possibly depending on the tree considered and further mathematical objects. These dependencies will be specified 
%by using a superscript as well as 
by writing $\alpha(...)$ and listing them in the brackets.  We shall use $\beta_\epsilon$ to denote constants, which converge as $\epsilon\to 0$. We shall again 
%use superscripts and 
write $\beta_\epsilon(...)$ in order to specify dependencies. Importantly, $\alpha$ and $\beta_\epsilon$ will always be constants and not functions!

%\begin{remark}\label{rem:why full space noise}
We can and shall choose renormalisation functions as if we were working with white noise in the full space instead of periodic noise.
One can check that 
this only affects the constants $\beta_\epsilon$ in the examples that follow.\footnote{This is consistent with the idea that `renormalisation is local'.}
%\end{remark}

\subsection{Heat Kernel regularisation}\label{subsec:heat kernel reg}
In this subsection,
when working with spatial white noise $\xi\in \mathcal{D}'(\mathbb{R}^d)$, we use the following approximation of the noise
\begin{equation}\label{eq:heat kernel regularisation spatial white noise}
\xi_\epsilon(x,t)= \int_{\mathbb{R}^d} \Gamma(x,t, \zeta, t-\epsilon^2) \xi(\zeta) \,d\zeta\ .
\end{equation}
Note in particular that $\xi_\epsilon\in \mathcal{C}^\infty(\mathbb{R}\times \mathbb{R}^d)$ is a function of space and time, despite the original noise being constant in time.

When working with space time white noise $\xi\in \mathcal{D}'( \mathbb{R}^{d+1})$ we fix an additional smooth compactly supported even function $\phi\in \mathcal{C}^\infty(\mathbb{R})$ such that $\int \phi =1$ and set
\begin{equation}\label{eq:heat kernel regularisation space time white noise}
\xi_\epsilon(x,t)= \int_{\mathbb{R}} \phi^\epsilon(t-\tau) \int_{\mathbb{R}^d}\Gamma(x,t, \zeta, t-\epsilon^2) \xi(\zeta,\tau) \,d\zeta d\tau\ ,
\end{equation}
where $\phi^\epsilon(t) :=  \frac{1}{\epsilon^2}\phi \left( \frac{t}{\epsilon^2} \right) $. 

\begin{remark}
Note that the exponent $2$ of $\epsilon^2$ in \eqref{eq:heat kernel regularisation spatial white noise} and \eqref{eq:heat kernel regularisation spatial white noise} is included so that the regularisation respects parabolic scaling in the sense of \cite{Hai14}.
\end{remark}
\subsubsection{Counterterms for g-PAM}\label{sec:PAM}
%Thus $$\bfPi \Xi (x,t)= \int_0^t \int_{\mathbb{R}^d} \Gamma(x,t, \zeta, \tau)\kappa(t-\tau) \xi (\zeta) d\zeta d\tau $$ 
Here we consider $d=2$ (but since a large part of the computations is generic we shall often still use the letter $d$). 
Working with $\xi_\epsilon$ as in \eqref{eq:heat kernel regularisation spatial white noise}, one finds that
\begin{align*}
\bfPi^{\epsilon} \<1> (x,t) &= \int_{\mathbb{R}} \int_{\mathbb{R}^d} \Gamma(x,t, \zeta, \tau) \kappa(t-\tau)\xi_\epsilon(\zeta,\tau) d\zeta d\tau \\
&=\int_{\mathbb{R}} \int_{\mathbb{R}^d}\int_{\mathbb{R}^d} \Gamma(x,t; \zeta, \tau) 
 \Gamma(\zeta,\tau, \eta, \tau-\epsilon^2) \xi(\eta) \,d\zeta d\eta \  \kappa(t-\tau) d\tau \\
 &=\int_{\mathbb{R}} \int_{\mathbb{R}^d} \Gamma(x,t ; \eta, \tau-\epsilon^2) \xi(\eta) \, d\eta \   \kappa(t-\tau) d\tau \\
 &=\int_{\mathbb{R}} \int_{\mathbb{R}^d} \Gamma^*( \eta, \tau-\epsilon^2; x,t ) \xi(\eta) \, d\eta  \ \kappa(t-\tau) d\tau 
\end{align*}
and therefore
\begin{align*}
\E\left[\bfPi^{\epsilon} \<Xi2> (x,t)\right]
%&= \int_{\mathbb{R}^d} \Gamma(x,t;\eta,t-\epsilon^2)\int \Gamma(x,t;\eta,\tau-\epsilon^2) d\eta\  \kappa(t-\tau) d\tau \\
&=\int_{\mathbb{R}^d} \Gamma^*(\eta,t-\epsilon^2;x,t)\int_{\mathbb{R}} \Gamma^*(\eta,\tau-\epsilon^2;x,t) d\eta\  \kappa(t-\tau) d\tau\ .
\end{align*}
After writing $\Gamma^*= Z_0^* + Z_{\geq 1}^*$ we investigate the contribution of each term to the above expected value. By Corollary~\ref{cor:bound_on_Z_>nu} the contributions involving $Z_{\geq 1}^*$ converge as $\epsilon\to 0$, while
the contribution only involving $ Z_0^*$ is given by 
 \begin{align*}
&\int_{\mathbb{R}} \int_{\mathbb{R}^d}  Z_0^*(\eta,t-\epsilon^2;x,t) Z_0^*(\eta,\tau-\epsilon^2;x,t) d\eta\  \kappa(t-\tau) d\tau\\
 &= 
\int_{\mathbb{R}} \int_{\mathbb{R}^d} 
\frac{
 C(x,t)}{(\epsilon^2)^{d/2}} \exp\left(- \frac{\vartheta^{(x,t)}(\eta-x) }{4\epsilon^2} \right) 
\frac{
 C(x, t)}{ (t-\tau+\epsilon^2)^{d/2}}
  \exp\left(- \frac{\vartheta^{(x,t)}(\eta-x) }{4(t-\tau+\epsilon^2)} \right) d\eta\  \kappa(t-\tau) d\tau\\
&= 
\int_{\mathbb{R}} \int_{\mathbb{R}^d} 
 C(x,t)^2(\epsilon^2)^{-d/2}(t-\tau+\epsilon^2)^{-d/2}
  \exp\left(- \frac{\vartheta^{(x,t)}(\eta-x) }{4\epsilon^2} 
 - \frac{\vartheta^{(x,t)}(\eta-x) }{4(t-\tau+\epsilon^2)} \right)   d\eta\  \kappa(t-\tau) d\tau \\
 &= 
\int_{\mathbb{R}} \int_{\mathbb{R}^d} 
 C(x,t)^2(\epsilon^2)^{-d/2}(t-\tau+\epsilon^2)^{-d/2} \exp\left(
 - \frac{\vartheta^{(x,t)}(\eta-x) }{4\frac{(t-\tau+\epsilon^2)\epsilon}{t-\tau +2\epsilon^2}} \right)  d\eta\  \kappa(t-\tau) d\tau  \\
 &= 
\int_{\mathbb{R}}
 C(x,t) (\epsilon^2)^{-d/2}(t-\tau+\epsilon^2)^{-d/2} \left( \frac{(t-\tau+\epsilon)\epsilon^2}{t-\tau +2\epsilon^2} \right)^{d/2}   \kappa(t-\tau) d\tau\\
% &= 
%\int_{\mathbb{R}}
% C(x,t) ({t-\tau +2\epsilon^2} )^{-d/2} \kappa(t-\tau) d\tau\\
 &=  C(x,t) 
\int_{\mathbb{R}}
({\tau +2\epsilon^2} )^{-d/2} \kappa(\tau) d\tau\\
 &=C(x,t)\left(2|\log (\epsilon)| + \beta_\epsilon(\kappa) \right) \ .
 \end{align*}
Next, we turn to the symbols $\{\<b2>_{i,j}\}_{i,j=1}^2$, where $\<bI>_i$ stands for the abstract lift of the derivative of the kernel $\<I>$ in the $i$-th direction.
Thus, we find using \eqref{eq:composition}
$$\bfPi^{\epsilon} \<b1>_{i}(x,t)=\int_{\mathbb{R}} \int_{\mathbb{R}^d} \partial_{x_i} \Gamma^*( \eta, \tau-\epsilon^2; x,t ) \xi(\eta) \, d\eta \ \kappa(t-\tau) d\tau \, $$
which implies that 
\begin{align*}
\E\left[\bfPi^{\epsilon} \<b2>_{ij} (x,t)\right] &=\int_{\mathbb{R}}\int_{\mathbb{R}} \int_{\mathbb{R}^d}\big( \partial_{x_i} \Gamma^*( \eta, \tau-\epsilon^2; x,t )  \big)  \big(  \partial_{x_j} \Gamma^*( \eta, \sigma-\epsilon^2; x,t )  \big)  d\eta   \kappa(t-\tau) d\tau \kappa(t-\sigma) d\sigma
\end{align*}
Writing $\partial_{x_i}\Gamma^*= \partial_{x_i}Z_0^* + \partial_{x_i}Z_{\geq 1}^*$ it follows from 
Remark~\ref{rem:upper bound components derivatives} that any contribution involving $Z_{\geq 1}^*$ to $\E\left[\bfPi^{\epsilon} \<b2>_{ij} (x,t)\right] $ converges to a finite limit. Finally, by Lemma~\ref{lem:bound on Z_{0;1}} the only divergent contribution to $\E\left[\bfPi^{\epsilon} \<b2>_{ij}  (x,t)\right]$ is 
\begin{align*}
&\int_{\mathbb{R}}\int_{\mathbb{R}} \int_{\mathbb{R}^d} Z_{0;i}^*( \eta, \tau-\epsilon^2; x,t )   Z_{0;j}^*( \eta, \sigma-\epsilon^2; x,t ) \, d\eta \  \kappa(t-\tau) d\tau\ \kappa(t-\sigma) d\sigma \\
&= \sum_{l,k} a^{i,l}a^{j,k} \int_{\mathbb{R}}\int_{\mathbb{R}} \left( \int_{\mathbb{R}^d} \frac{(\eta-x)_l(\eta-x)_k }{4(t-\tau +\epsilon^2)(t-\sigma +\epsilon^2)}  Z_0^*( \eta, \tau-\epsilon^2; x,t ) Z_{0}^*( \eta, \sigma-\epsilon^2; x,t ) \, d\eta \right) \\
&\qquad \qquad\qquad \qquad\qquad \times  \kappa(t-\tau) d\tau\ \kappa(t-\sigma) d\sigma \ .
\end{align*}
We calculate 
\begin{align*}
&\int_{\mathbb{R}^d} (\eta-x)_l(\eta-x)_k  Z_0^*( \eta, \tau-\epsilon^2; x,t ) Z_{0}^*( \eta, \sigma-\epsilon^2; x,t ) \, d\eta\\
%&=    \int_{\mathbb{R}^d} C(x, t)^2 \frac{(\eta-x)_l(\eta-x)_k} { (t-\tau+\epsilon^2)^{d/2}(t-\sigma+\epsilon^2)^{d/2}}
%  \exp\left(- \frac{\vartheta^{(x,t)}(\eta-x) }{4(t-\tau+\epsilon^2)}  - \frac{\vartheta^{(x,t)}(\eta-x) }{4(t-\sigma+\epsilon^2)} \right)d\eta\\
  &=    \int_{\mathbb{R}^d} C(x, t)^2 \frac{(\eta-x)_l(\eta-x)_k} { (t-\tau+\epsilon^2)^{d/2}(t-\sigma+\epsilon^2)^{d/2}}
  \exp\left(- \frac{\vartheta^{(x,t)}(\eta-x) }{4\frac{(t-\tau+\epsilon^2)(t-\sigma+\epsilon^2)}{2(t+\epsilon^2) -\sigma-\tau}}  \right) d\eta\\
%&=    \int_{\mathbb{R}^d} C(x, t)^2 \frac{\eta_l\eta_k} { (t-\tau+\epsilon^2)^{d/2}(t-\sigma+\epsilon^2)^{d/2}}
%  \exp\left(- \frac{\vartheta^{(x,t)}(\eta) }{4\frac{(t-\tau+\epsilon^2)(t-\sigma+\epsilon^2)}{2(t+\epsilon^2) -\sigma-\tau}}  \right) d\eta\\  
 &= \left(\frac{(t-\tau+\epsilon^2)(t-\sigma+\epsilon^2)}{2(t+\epsilon^2) -\sigma-\tau}\right)^{d/2+1} \int_{\mathbb{R}^d} C(x, t)^2 \frac{\eta_l\eta_k} { (t-\tau+\epsilon^2)^{d/2}(t-\sigma+\epsilon^2)^{d/2}}
  \exp\left(- \frac{\vartheta^{(x,t)}(\eta) }{4}  \right) d\eta\\    
   &=    \frac{(t-\tau+\epsilon^2)(t-\sigma+\epsilon^2) C(x, t)^2}{(2(t+\epsilon^2) -\sigma-\tau)^{d/2+1}}
  \int_{\mathbb{R}^d} \eta_l \eta_k\exp\left(- \frac{\vartheta^{(x,t)}(\eta) }{4}  \right) d\eta\\
  &=   2  {a_{l,k}(x,t)} C(x, t)\frac{ (t-\tau+\epsilon^2)(t-\sigma+\epsilon^2)}{(2(t+\epsilon^2) -\sigma-\tau)^{d/2+1}}
\end{align*}
and finally find
\begin{align*}
&\int_{\mathbb{R}}\int_{\mathbb{R}} \int_{\mathbb{R}^d} Z_{0;i}^*( \eta, \tau-\epsilon^2; x,t )   Z_{0;j}^*( \eta, \sigma-\epsilon^2; x,t ) \, d\eta \ \kappa(t-\tau) d\tau\ \kappa(t-\sigma) d\sigma \\
&= \sum_{l,k} a^{i,l}a^{j,k} \int_{\mathbb{R}}\int_{\mathbb{R}} \frac{1}{4}  \frac{ 2a_{l,k} C(x, t)}{(2(t+\epsilon^2) -\sigma-\tau)^{d/2+1}}
  \ \kappa(t-\tau) d\tau\ \kappa(t-\sigma) d\sigma \\
  &= \frac{a^{i,j} C(x, t)}{2} \int_{\mathbb{R}}\int_{\mathbb{R}} {(2\epsilon^2 +\sigma+\tau)^{-d/2-1}} \ \kappa(\tau) d\tau\ \kappa(\sigma) d\sigma\\
   &= a^{i,j}(x,t) C(x, t) \left(|\log(\epsilon)| + \beta_\epsilon( \kappa) \right) \ .
\end{align*}
%\red{Thus we obtain that contribution of $Z_{0; i}$ to $E[\bfPi^{\epsilon} (DI\Xi\otimes DI\Xi)(x,t)]$ is given by
% $$ C(x, t)A^{-1}(x,t)\left(a \log(\epsilon) +b \right)$$ with $b$ depending on the chosen time cut-off, since all terms involving $Z_{\geq 1}$ or $(\partial_i Z_0 -Z_{0;i} )$ are seen to be convergent.
%}

%\gray{
%since
%$$\partial_i Z_0^*( \eta, \tau-\epsilon; x,t )= -\frac{\sum_l (a^{i,l}+ a^{l,i})(\eta-x)_l +\langle(\eta-x), (D_xA(x,t))(\eta-x)\rangle }{4(t-\tau +\epsilon)} Z_0^*( \eta, \tau-\epsilon; x,t ) $$
% 
% $$Z^*(x,t;\zeta, \tau)= C(\zeta, \tau) (\tau-t)^{-d/2} \exp\left( \frac{\vartheta^{(\zeta,\tau)}(x-\zeta) }{4(\tau-t)} \right) $$
%
%Using \eqref{eq:bound_on_Z_>nu}
%we note that the contribution only involving $Z_{\geq 1}^*$ stays bounded in $\epsilon$. 
%Let us investigate the contribution from  $Z_0^*$ which is given by 
%$$\int_{\mathbb{R}^d} \left(\int_0^t Z_0^*( \eta, \tau-\epsilon; x,t ) \mathbf{1}_{\tau>0} d\tau\right)  
%\left(\int_0^t Z_0^*( \eta, \sigma-\epsilon; x,t )\mathbf{1}_{\sigma>0} d\sigma \right) \xi(\eta) \, d\eta  $$
%
%
%\begin{align*}
%E\left[(\bfPi^{\epsilon} I\Xi (x,t)^2\right] &= \int_{\mathbb{R}^d} \left(\int_0^t \Gamma^*( \eta, \tau-\epsilon; x,t ) \mathbf{1}_{\tau>0} d\tau\right)  
%\left(\int_0^t \Gamma^*( \eta, \sigma-\epsilon; x,t )\mathbf{1}_{\sigma>0} d\sigma \right) \xi(\eta) \, d\eta  
%\end{align*}
%
%Next we turn to 
%
%The relevant trees are $\Xi(I\Xi)$ and $I_i(\Xi) I_j(\Xi)$, we have for $\star\in \mathbb{R}^{d+1}$
%\begin{align*}
%E[\Pi^{(\eps)}_\star \Xi(I\Xi)] = E[K(\star, \bar{z}) \rho(\bar)]...
%\end{align*}
%}

\subsubsection{Counterterms for $\Phi^4_2$}\label{sec:phi42}
%This time we work with space time white noise $\bfPi\Xi(t,x)= \xi(x,t)$ and thus an additional regularisation in the time direction given by 
% $\phi^\epsilon(t) :=  \frac{1}{\epsilon^2}\phi \left( \frac{t}{\epsilon^2} \right) $, so
%we shall take the following approximation of the noise
%$$\xi_\epsilon(x,t)= \int_{\mathbb{R}} \phi^\epsilon(t-\tau) \int_{\mathbb{R}^d}\Gamma(x,t, \zeta, t-\epsilon^2) \xi(\zeta,\tau) \,d\zeta d\tau\ .$$
This time we work with space time white noise and $\xi_\epsilon$ as in \eqref{eq:heat kernel regularisation space time white noise} and again $d=2$. We find that
\begin{align*}
\bfPi^{\epsilon} \<1> (x,t) &= \int_\mathbb{R} \int_{\mathbb{R}^d} \Gamma(x,t, \zeta, \tau) \xi_\epsilon(\zeta,\tau) d\zeta \ \kappa(t-\tau) d\tau \\
%&=\int_\mathbb{R}\int_\mathbb{R} \int_{\mathbb{R}^d}\int_{\mathbb{R}^d} \Gamma(x,t; \zeta, \tau) 
% \phi^\epsilon(\tau-\sigma) \Gamma(\zeta,\tau, \eta, \tau-\epsilon^2) \xi(\eta,\sigma) \,d\eta d\sigma d\zeta \ \kappa(t-\tau) d\tau\\
&=\int_\mathbb{R}\int_\mathbb{R} \int_{\mathbb{R}^d} \Gamma(x,t; \eta, \tau-\epsilon^2) 
 \phi^\epsilon(\tau-\sigma) \xi(\eta,\sigma) \,d\eta d\sigma \ \kappa(t-\tau) d\tau\ .
\end{align*}
Thus, writing $\left(\phi^\epsilon\right)^{*2} ( \tau-\tau' ):= \int_\mathbb{R}
\phi^\epsilon(\tau-\sigma) \phi^\epsilon(\tau '-\sigma) $ find
\begin{align}
&\E [\bfPi^{\epsilon} \<2> (x,t)] \label{eq:cherry general form}\\
%&= \int_\mathbb{R}\int_\mathbb{R}\int_\mathbb{R}
%\phi^\epsilon(\tau-\sigma) \phi^\epsilon(\tau '-\sigma) 
% \int_{\mathbb{R}^d} \Gamma(x,t; \eta, \tau-\epsilon^2) 
% \Gamma(x,t; \eta, \tau'-\epsilon^2) 
%\,d\eta d\sigma  \ \kappa(t-\tau) d\tau \ \kappa(t-\tau') d\tau' \nonumber\\
&= \int_\mathbb{R}\int_\mathbb{R}
\left(\phi^\epsilon\right)^{*2} ( \tau-\tau' )
 \int_{\mathbb{R}^d} \Gamma^*(\eta, \tau-\epsilon^2;x,t) 
 \Gamma^*( \eta, \tau'-\epsilon^2;x,t; ) 
\,d\eta \ \kappa(t-\tau) d\tau \ \kappa(t-\tau') d\tau' \ . \nonumber 
\end{align}
Writing $\Gamma^*= Z^*_0+ Z^*_{\geq 1}$, again all contributions involving 
$Z^*_{\geq 1}$ converge to a finite limit as $\epsilon \to 0$,
while the term only involving $Z_0^*$ is given by
\begin{align*}
 &\int_\mathbb{R}\int_\mathbb{R}
\left(\phi^\epsilon\right)^{*2} ( \tau-\tau' )
 \int_{\mathbb{R}^d} Z_0^*(\eta, \tau-\epsilon^2;x,t) 
 Z_0^*( \eta, \tau'-\epsilon^2;x,t; ) 
\,d\eta
 \ \kappa(t-\tau) d\tau \ \kappa(t-\tau') d\tau'\\
&=\int_\mathbb{R}\int_\mathbb{R}
\left(\phi^\epsilon\right)^{*2} ( \tau-\tau' )\kappa(t-\tau)  \kappa(t-\tau')\\
&\quad\times \left(
\int_{\mathbb{R}^d} 
\frac{
C(x,t)^2}{ (t-\tau+\epsilon^2)^{d/2}(t-\tau'+\epsilon^2)^{d/2}}
\exp\left( - \frac{\vartheta^{(x,t)}(x-\eta)}{4 (t-\tau+\epsilon^2)} - \frac{\vartheta^{(x,t)}(x-\eta)}{4 (t-\tau'+\epsilon^2)} \right)
\,d\eta  \right)  \  d\tau d\tau'\\
 &=C(x,t) \int_\mathbb{R}\int_\mathbb{R}
\left(\phi^\epsilon\right)^{*2} ( \tau-\tau' ) 
 (2(t+\epsilon^2) -\tau -\tau')^{-d/2}
   \kappa(t-\tau) d\tau \ \kappa(t-\tau') d\tau'\\
    &=C(x,t) \int_\mathbb{R}\int_\mathbb{R}
\left(\phi^\epsilon\right)^{*2} ( \tau-\tau' ) 
 (2\epsilon^2+ \tau +\tau')^{-d/2}
   \kappa(\tau) d\tau \ \kappa(\tau') d\tau'\\
 &= C(x,t) \left(|\log(\epsilon)|  + \beta_\epsilon(\kappa, \phi)\right) \ ,
\end{align*}
where the last equality can be seen by substitution $(\kappa, \kappa') \mapsto (\epsilon^2 \kappa,\epsilon^2 \kappa')$.

\subsubsection{Counterterms for $\Phi^4_3$}\label{sec:phi43}
Here $d=3$. Then, for $
\E [\bfPi^{\epsilon} \<2> (x,t)]$
one obtains the identity \eqref{eq:cherry general form}, exactly as above.
This time, we write
$$\Gamma^* = Z_0^* + \bar{Z}_1^* + R_1^* + Z^*_{\geq 2} \, $$
where we used Lemma~\ref{lem:reflection} to write $Z_1^*=\bar{Z}_1^* + R_1^*$. In order to identify the divergent contributions to $\E [\bfPi^{\epsilon} \<2> (x,t)]$ we observe the following.
\begin{itemize}
\item The contributions of the terms involving either $R_1^*$ or $Z^*_{\geq 2}$ converge to a finite limit as $\epsilon\to 0$.
\item Similarly, the contribution only involving $\bar{Z}^*_1$ converges to a finite limit.
\item The contribution of the cross term between $Z_0^*$ and $\bar{Z}^*_1$ given by
\begin{equation*}
\int_\mathbb{R}\int_\mathbb{R}\int_\mathbb{R}
\phi^\epsilon (\tau-\sigma ) 
\phi^\epsilon ( \tau'-\sigma)
 \int_{\mathbb{R}^d} Z_0^*(\eta, \tau-\epsilon;x,t) 
 \bar{Z}_1^*( \eta, \tau'-\epsilon;x,t) 
\,d\eta d\sigma  d\tau d\tau'  \ .
\end{equation*}
vanishes due to \eqref{eq:reflection_symmetry}.
\end{itemize}
Thus, the only remaining contribution to $\E [\bfPi^{\epsilon} (I\Xi)^2 (x,t)]$ is as above given by 
\begin{align*}
&\int_\mathbb{R}\int_\mathbb{R}\int_\mathbb{R}
\phi^\epsilon ( \tau-\sigma ) 
\phi^\epsilon ( \tau'-\sigma )
 \int_{\mathbb{R}^d} Z_0^*(\eta, \tau-\epsilon^2;x,t) 
 Z_0^*( \eta, \tau'-\epsilon^2;x,t; ) 
\,d\eta
 d\sigma  \ \kappa(t-\tau) d\tau \ \kappa(t-\tau') d\tau'\\
 &=C(x,t) \int_\mathbb{R}\int_\mathbb{R}
\left(\phi^\epsilon\right)^{*2} ( \tau-\tau' ) 
 (2\epsilon^2 +\tau +\tau')^{-d/2}\ 
   \kappa(\tau) d\tau \ \kappa(\tau') d\tau'\\
 &= C(x,t) \left( \frac{\alpha( \phi)}{\epsilon} +  \beta_\epsilon (\kappa, \phi)\right)\ .
\end{align*}
For the $\phi^4_3$-equation there is an additional divergence coming from the tree $\<22>$, see \cite{Hai14}. It is given by
\begin{align*}
&\int_\mathbb{R} \int_{\mathbb{R}^d} \Gamma(x,t; \bar{x},\bar{t}) \left(\tilde{\Gamma}_\epsilon(x,t, \bar{x},\bar{t})\right)^2 \ d\bar{x} \kappa(t-\bar{t}) d\bar{t}\\
%%&= \int_\mathbb{R} \int_{\mathbb{R}^d} \Gamma^*(\bar{x},\bar{t};x,t, ) \left(\tilde{\Gamma}_\epsilon(x,t, \bar{x},\bar{t})\right)^2 \ d\bar{x\kappa(t-\bar{t}) d\bar{t}\\
  &= \int_\mathbb{R} \int_{\mathbb{R}^d} \Gamma^*(\bar{x}+x,\bar{t}+t;x,t ) \left(\tilde{\Gamma}_\epsilon(x,t, \bar{x}+x,\bar{t}+t)\right)^2 \ d\bar{x} \kappa(t-\bar{t}) d\bar{t}\
\end{align*}
where
\begin{align*}
\tilde{\Gamma}_\epsilon(x,t, \bar{x},\bar{t})&= \int_\mathbb{R}\int_\mathbb{R}
(\phi^\epsilon)^{*2} (\tau-\tau' ) 
 \int_{\mathbb{R}^d} \Gamma (x,t;\eta, \tau-\epsilon^2) 
 \Gamma ( \bar{x},\bar{t};\eta, \tau'-\epsilon^2) 
\,d\eta \ \kappa(t-\tau) d\tau \ \kappa(\bar{t}-\tau') d\tau'\\
&= \int_\mathbb{R}\int_\mathbb{R}
(\phi^\epsilon)^{*2} (\tau-\tau' ) 
 \int_{\mathbb{R}^d} \Gamma^* (\eta, \tau-\epsilon^2;x,t) 
 \Gamma^* (\eta, \tau'-\epsilon^2; \bar{x},\bar{t}) 
\,d\eta \ \kappa(t-\tau) d\tau \ \kappa(\bar{t}-\tau') d\tau'\ .
\end{align*}
Writing $\Gamma^*= Z_0^* +  Z_{\geq 1}^*$ one finds that the only divergent summand is given by
$$\int_\mathbb{R} \int_{\mathbb{R}^d} Z^*(\bar{x},\bar{t};x,t ) \left(\tilde{Z}^*_\epsilon(\bar{x},\bar{t};x,t )\right)^2 \ d\bar{x} \kappa(t-\bar{t}) d\bar{t}$$
where 
\begin{align*}
&\tilde{Z}^*_\epsilon(x,t, \bar{x},\bar{t})\\
&:= \int_\mathbb{R}\int_\mathbb{R}
(\phi^\epsilon)^{*2} (\tau-\tau' ) 
 \int_{\mathbb{R}^d} Z^{*} (\eta, \tau-\epsilon^2;x,t;) 
 Z^* ( \eta, \tau'-\epsilon^2;\bar{x},\bar{t}) 
\,d\eta \ \kappa(t-\tau) d\tau \ \kappa(\bar{t}-\tau') d\tau'\\ 
&= \int_\mathbb{R}\int_\mathbb{R}
(\phi^\epsilon)^{*2} (\tau-\tau' ) 
 \int_{\mathbb{R}^d} Z^{*} (\eta, \tau-\epsilon^2;x,t;) 
 C(\bar{x},\bar{t}) w^{(\bar{x},\bar{t})} ( \eta, \tau'-\epsilon^2;\bar{x},\bar{t}) 
\,d\eta \ \kappa(t-\tau) d\tau \ \kappa(\bar{t}-\tau') d\tau' \ .
\end{align*}
\begin{lemma}\label{lem:heat kernel recentering}
It holds that
 $$|C(\bar{x},\bar{t}) w^{(\bar{x},\bar{t})} ( \eta, \tau;\bar{x},\bar{t}) - 
 C({x},{t}) w^{({x},{t})} ( \eta, \tau;\bar{x},\bar{t}) |\lesssim
 (|x-\bar{x}|+ |t-\bar{t}|^{1/2}) (t-\tau)^{-d/2}  \exp\left( -\frac{\lambda^*_0 |\bar{x}-\eta|^2}{4(\bar{t}-\tau)} \right),$$
 where the implicit constant depends on $A$.
 \end{lemma}
Inspired by the above lemma, set
\begin{align*}
\hat{Z}^*_\epsilon(x,t, \bar{x},\bar{t})
&:=\int_\mathbb{R}\int_\mathbb{R}
(\phi^\epsilon)^{*2} (\tau-\tau' ) \kappa(t-\tau)\kappa(\bar{t}-\tau') \\
&\qquad \times
\left( \int_{\mathbb{R}^d} Z^{*} (\eta, \tau-\epsilon^2;x,t) 
 C({x},{t}) w^{({x},{t})} ( \eta, \tau'-\epsilon;\bar{x},\bar{t}) 
\,d\eta \right) \  d\tau  d\tau'\\
&=\int_\mathbb{R}\int_\mathbb{R}
(\phi^\epsilon)^{*2} (\tau-\tau' )  \kappa(t-\tau)\kappa(\bar{t}-\tau')\\
& \qquad \times \left(
 \int_{\mathbb{R}^d} C({x},{t})^2 w^{({x},{t})}(x,\tau-\epsilon^2;\eta, t) 
 w^{({x},{t})} ( \eta, \tau'-\epsilon^2;\bar{x},\bar{t}) 
\,d\eta \right)\  d\tau  d\tau'\\
&=\int_\mathbb{R}\int_\mathbb{R}
(\phi^\epsilon)^{*2} (\tau-\tau' ) \kappa(t-\tau)\kappa(\bar{t}-\tau')\\
&\qquad \times \left(  C({x},{t}) w^{({x},{t})}(x,\tau+\tau'-2\epsilon^2;\bar{x}, t+\bar{t}) 
\,  \right)\  d\tau  d\tau'\\
&=\int_\mathbb{R}\int_\mathbb{R}
 (\phi^\epsilon)^{*2} (\tau-\tau' ) \kappa(t-\tau)\kappa(\bar{t}-\tau')
Z^*(\bar{x},\tau+\tau'-2\epsilon^2-\bar{t};{x}, t) 
\, d\tau  d\tau'
\end{align*}
and observe that the remaining diverging summand is given by 
\begin{align*}
\int_\mathbb{R} \int_{\mathbb{R}^d} Z^*(\bar{x},\bar{t};x,t ) \left(\hat{Z}^*_\epsilon(\bar{x},\bar{t};x,t )\right)^2 \ d\bar{x} \kappa(t-\bar{t}) d\bar{t}
&=C(x,t)^2 \left(\alpha |\log(\epsilon)| + \beta_\epsilon (\phi, \kappa)\right)\ .
\end{align*}

\begin{remark}\label{rem:phi34cherry}
Let us observe at this point that in $d\geq 4$ dimensions the identification of the divergent contributions to \eqref{eq:cherry general form}
would require further analysis. For $d=4$ a decomposition of the form
$$\Gamma^* = Z_0^* + \bar{Z}_1^* + R_1^* + \bar{Z}^*_{ 2} + R_2^* + \bar{Z}^*_{\geq 3}\ . $$
 then yields additional possible divergences, which
%\begin{itemize}
%\item only involving $\bar{Z}_1^*$
%\item from the cross term involving $Z_0^*$ and $R_1^*$ 
%\item from the cross term of $Z_0^*$ and $\bar{Z}^*_{ 2}$ .
%\end{itemize}
 are explicit local expressions involving $a$, derivatives thereof, as well as possibly $b$ and $c$. Despite the contributions involving $b$ and $c$ not being expected to diverge, since they can be absorbed into the right hand side of the equation, 
one can see that an additional divergence does appear for the non-translation invariant $\phi^{3}_4$-equation. This follows from Remark~\ref{rem:riemannian} and the fact that for the geometric $\phi^{3}_4$-equation in \cite{HS23m} there is an additional logarithmic divergence proportional to the scalar curvature. Thus, in particular, the scaling heuristic of Section~\ref{sec:scaling heuristic} does not apply to this equation.
\end{remark}

\subsubsection{Counterterms for KPZ}\label{sec:KPZ}
For the KPZ-equation, where $d=1$ we first write 
\begin{equation}\label{eq:kernel_dec_kpzs}
\partial_1 \Gamma^*= Z_{0;1}^* + \partial_1 \bar{Z}_{ 1}^*   +R^*_{0,1} + \partial R^*_{1} + \partial_1 Z_{\geq 2}^* ,
\end{equation}
where we used the decompositions from Lemma~\ref{lem:bound on Z_{0;1}} in the case $d=1$ and Lemma~\ref{lem:reflection}. 
After adopting the usual convention for the KPZ-equation that 
$\<I>$ stands for the derivative of the heat kernel,
we make the following observation about contribution of each term in \eqref{eq:kernel_dec_kpzs} to $\E\left[\bfPi^{\epsilon} \<2> (x,t)\right]$:
\begin{enumerate}
\item Any term involving either $\partial R^*_{1}$ or $\partial_1 Z_{\geq 2}^*$ converges as $\epsilon\to 0$.
\item The cross term involving $\partial_1 \bar{Z}_{ 1}^* $  and $R^*_{0,1}$ converges as $\epsilon\to 0$.
\item The cross terms involving $Z_{0;1}^*$ and either $\partial_1 \bar{Z}_{ 1}^* $ or $R^*_{0,1}$ vanish, since $Z_{0;1}^*$ is odd with respect to $\mathfrak{R}$ defined in \eqref{eq:reflection}, while the latter two are even.
\end{enumerate}
Thus, the only remaining term involves just $Z_{0;1}^*$ and is given by
\begin{align*}
 &\int_\mathbb{R}\int_\mathbb{R}
(\phi^\epsilon)^{*2} ( \tau-\tau' ) 
 \int_{\mathbb{R}^d} Z_{0;1}^*(\eta, \tau-\epsilon^2;x,t) 
 Z_{0;1}^*( \eta, \tau'-\epsilon^2;x,t; ) 
\,d\eta
 d\sigma  \ \kappa(t-\tau) d\tau \ \kappa(t-\tau') d\tau'\\
& = a^{-1}(x,t) C(x, t)\left( \frac{\alpha'(\phi)}{\epsilon} + \beta'_\epsilon( \phi, \kappa) \right) \\
& = a^{-3/2}(x,t) \left( \frac{\alpha(\phi)}{\epsilon} + \beta_\epsilon( \phi, \kappa) \right) \ . 
\end{align*}

For the remaining counterterms appearing due to the diagrams  $\<11>,\  \<22j>,  \<211>$
 we observe the following:
\begin{enumerate}
\item $\E[\bfPi^\epsilon \<11>(p,t)]$ converges since the only possibly divergent term vanishes due to reflection properties. (This term ``would'' come with the function $a^{-3/2}(x,t) C(x,t)\sim a^{-2}(x,t)$ in front of the divergence.)
\item Both logarithmic divergences corresponding to $\<22j>,  \<211>$ have a pre-factor proportional to $C^2(x,t) a^{-3}(x,t)\sim a^{-4}(x,t)$.
\end{enumerate}
\begin{remark}
One can show that, as for the usual KPZ equation, the logarithmic divergences due to $\<22j>,  \<211>$ cancel, c.f.\ \cite[Theorem~6.5]{HQ18} and \cite[Lemma~6.5]{Hai13}.
\end{remark}

\subsection{More general `covariant' regularisations}\label{sec:other regularisation}

In this section we show that using the heat kernel $\Gamma$ in order to regularize the noise in the spatial direction, while convenient, is not crucial;
its main advantage being that it is more robust, c.f.\ Remark~\ref{rem:fail}, and allows for simple explicit computations.
Here we consider a type of regularisation, which is `sufficiently covariant' for the equations considered,\footnote{We shall frequently use the expression`the equations considered' meaning the g-PAM, $\phi^4_2,$ $\phi^4_3$ and KPZ-equation.} allowing one to choose counterterms of the same form as when using heat kernel regularisation and in the limit obtain the same solutions.

Let $\rho\in \cC_c^\infty(\R_+)$ non-negative with support in $[0,1)$, all odd derivatives vanishing at the origin, and such that $\int_{\mathbb{R}^d} \rho (|x|^2) dx=1$. Set $$\rho^{(z, \epsilon)}(x)=\frac{1}{\epsilon^d \det(A(z))^{1/2}} \rho \left(\frac{\vartheta^{z} (x)}{\epsilon^2}\right)\ .$$
We shall now study the case when the noise $\xi_\epsilon\in \mathcal{C}^\infty ( \mathbb{R}^{d+1})$ is obtained as in \eqref{eq:heat kernel regularisation spatial white noise}, resp.~\eqref{eq:heat kernel regularisation space time white noise} but with $\Gamma(x,t, \zeta, t-\epsilon^2)$ replaced by $\rho^{((x,t),\epsilon)} (x-\zeta)$. 
In the former case, this reads
\begin{equation}\label{eq:general_regular_spatial}
\xi_\epsilon(x,t)= \int_{\mathbb{R}^d} \rho^{((x,t),\epsilon)} (x-\zeta) \xi(d\zeta) \ ,
\end{equation}
but for notational convenience we shall only discuss the case of space time white noise from now on.\footnote{Only spatial white noise can be treated similarly.}
We set for $\phi^\epsilon$ as in \eqref{eq:heat kernel regularisation space time white noise} 
$$ \varrho^\epsilon (x,t; \zeta, \tau)  = \phi^\epsilon(t-\tau) \rho^{((x,t),\epsilon)} (x-\zeta)  \ ,$$
then the mollified noise is given by
\begin{equation}\label{eq:gen_regularisation}
\xi_\epsilon(z) = \int \varrho^\epsilon (z; z') \xi(z')\, dz'\ . 
\end{equation}
%Finally set $$ \varrho^\epsilon (x,t; \zeta)  =  \frac{1}{\epsilon}\rho^{(x,t)} \left( \frac{x-\zeta}{\epsilon}\right) \ .$$
For later use, define $\varrho^{z,\epsilon} (x,t; \zeta,\tau)  =  \phi^\epsilon(t-\tau) \rho^{(z,\epsilon)} (x-\zeta)  \ .$ 
%and shall also write 
%$\varrho^\epsilon (x,t; \zeta,s)$ and $\varrho^{z,\epsilon} (x,t; \zeta,s)$ with the understanding that the function constant in the $s$ variable.

%Finally, regularise spacial white noise $\xi$ by
%$$\xi_\epsilon(x,t)= \int_\mathbb{R}^d \varrho^\epsilon (x,t;\zeta) \xi(d\zeta)$$ and space time white noise by 
%$$\xi_\epsilon(x,t)= \int_{\mathbb{R}\times \mathbb{R}^d} \phi^{\epsilon}(s)\varrho^\epsilon (x,t;\zeta) \xi(d\zeta, ds)\ .$$
% \frac{1}{\epsilon^2 \det(A(x,t)^{1/2}} \rho \left(\frac{\vartheta^{x,t} (x-\zeta)}{\epsilon^2}\right)\ ,

All expected values we need to renormalise can be represented by a directed graph $\CCG= (\CCV, \CCE)$ as follows.
We are given two types of edges $\CCE=\CCE_+\cup \CCE_0$, the former representing the kernel $\Gamma^*$ (resp.\ $\partial_i \Gamma^*$ for g-Pam or KPZ) and the latter representing instances of $\varrho^\epsilon$. To each vertex $v\in \CCV$ we associate a variable in $z_v\in \mathbb{R}^{d+1}$. The set $\CCV$ contains one special vertex $\star$ and we shall integrate over all variables in $\CCV_0= \CCV\setminus \{\star\}$. Thus to each such graph $\CCG$ we can associate a corresponding integral expression
\begin{equation}\label{eq:general_integral_involving_Gammas}
 \int_{z\in (\mathbb{R}\times \mathbb{R}^d)^{\CCV_0}} \prod_{e\in \CCE_+} \Gamma^* (z_{e_+},z_{e_-})   \prod_{e\in \CCE_0}  \varrho^\epsilon  (z_{e_+},z_{e_-})  dz\ ,
\end{equation}
% \int_{x\in (\mathbb{R}\times \mathbb{R}^d)^{\CCV_0}} \prod_{e\in \CCE_+} \Gamma^* (z_{e_+},z_{e_-})   \prod_{e\in \CCE_0}  \varrho^\epsilon  (z_{e_+},z_{e_-})  dx\ ,
%$$
respectively 
\begin{equation}\label{eq:general_integral_involving_nablaGammas}
 \int_{z\in (\mathbb{R}\times \mathbb{R}^d)^{\CCV_0}} \bigotimes_{e\in \CCE_+} \nabla_{x_{e_-}}\Gamma^* (z_{e_+},z_{e_-})   \prod_{e\in \CCE_0}  \varrho^\epsilon  (z_{e_+},z_{e_-})  dz\ .
\end{equation}
We introduce the following graphical notation:
\begin{enumerate}
\item $\tikz[baseline=-3] \node [root] {};$ represents the special node $\star$, that is an instance of the variable $z_\star$,
\item $\tikz[baseline=-0.1cm] \draw[arrho] (0,0) to (1,0);$
represents $e\in \CCE_{0}$ directed from $e_-$ to $e_+$, i.e.\ an instance of $\varrho^\epsilon  (z_{e_+},z_{e_-})$,
\item $\tikz[baseline=-0.1cm] \draw[kernel] (0,0) to (1,0);$
represents a directed edge $e\in \CCE_{+}$ directed from $e_-$ to $e_+$, i.e.\ an instance of $\Gamma^* (z_{e_+},z_{e_-})$ resp. $\nabla_{x_{e_-}}\Gamma^* (z_{e_+},z_{e_-})$.
\item $\tikz[baseline=-3] \node [dot] {};$ represents a variable in $\mathbb{R}\times \mathbb{R}^d$ to be integrated over.
\end{enumerate}
Thus, all divergent expressions for the equations considered can be represented by one of following diagrams\footnote{The first two in the case of g-PAM, being interpreted with the appropriate modification.}
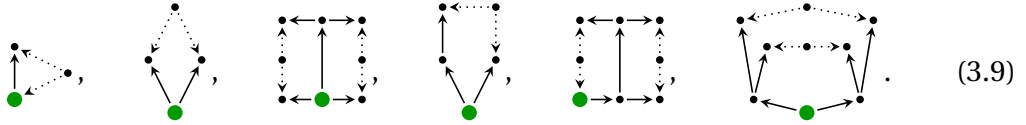
\begin{equation}\label{eq:graphs}
\begin{tikzpicture}[scale=0.35,baseline=0.6cm]
	\node at (-2,1)  [root] (left) {};
	\node at (-2,3)  [dot] (left1) {};
	\node at (0,2) [dot] (variable) {};
	
	\draw[kernel] (left) to  (left1);
	\draw[arrho] (variable) to (left1); 
	\draw[arrho] (variable) to (left); 
\end{tikzpicture}\;,
\qquad
\begin{tikzpicture}[scale=0.35,baseline=0.6cm]
	\node at (0,0.5)  [root] (root) {};
	\node at (-1,2.5)  [dot] (left) {};
	\node at (1,2.5)  [dot] (right) {};
	\node at (0,4.5) [dot] (variable) {};
	\draw[kernel] (root) to (left);
	\draw[kernel] (root) to (right);
	\draw[arrho] (variable) to (left); 
	\draw[arrho] (variable) to (right); 
\end{tikzpicture}\;,
\qquad
\begin{tikzpicture}[scale=0.35,baseline=0.6cm]
	\node at (0,1)  [root] (root) {};
	\node at (-1.5,1)  [dot] (left) {};
	\node at (1.5,1)  [dot] (right) {};
	\node at (-1.5,2.5)  [dot] (leftleft) {};
	\node at (1.5,2.5)  [dot] (rightright) {};
	\node at (-1.5,4)  [dot] (topleft) {};
	\node at (1.5,4)  [dot] (topright) {};
	\node at (0,4) [dot] (top) {};
	\draw[kernel] (root) to (left);
	\draw[kernel] (root) to (right);
	\draw[kernel] (root) to (top);
	\draw[kernel] (top) to (topleft); 
	\draw[kernel] (top) to (topright); 
	\draw[arrho]  (leftleft) to (topleft); 
	\draw[arrho]  (rightright) to (topright); 
	\draw[arrho] (leftleft) to (left) ; 
	\draw[arrho] (rightright) to (right) ; 
\end{tikzpicture}\;,
\qquad
\begin{tikzpicture}[scale=0.35,baseline=0.6cm]
	\node at (0,0.5)  [root] (root) {};
	\node at (-1,2.5)  [dot] (left) {};
	\node at (1,2.5)  [dot] (right) {};
	\node at (-1,4.5) [dot] (topleft) {};
	\node at (1,4.5) [dot] (topright) {};
	\draw[kernel] (root) to (left);
	\draw[kernel] (root) to (right);
	\draw[kernel] (left) to (topleft); 
	\draw[arrho] (topright) to (right); 
	\draw[arrho] (topright) to (topleft); 
\end{tikzpicture}\;,
\qquad
	\begin{tikzpicture}[scale=0.35,baseline=0.6cm]
	\node at (0,1)  [dot] (middle) {};
	\node at (-1.5,1)  [root] (left) {};
	\node at (1.5,1)  [dot] (right) {};
	\node at (-1.5,2.5)  [dot] (leftleft) {};
	\node at (1.5,2.5)  [dot] (rightright) {};
	\node at (-1.5,4)  [dot] (topleft) {};
	\node at (1.5,4)  [dot] (topright) {};
	\node at (0,4) [dot] (top) {};
	\draw[kernel] (left) to (middle);
	\draw[kernel] (middle) to (right);
	\draw[kernel] (middle) to (top);
	\draw[kernel] (top) to (topleft); 
	\draw[kernel] (top) to (topright); 
	\draw[arrho]  (leftleft) to (topleft); 
	\draw[arrho]  (rightright) to (topright); 
	\draw[arrho] (leftleft) to (left) ; 
	\draw[arrho] (rightright) to (right) ; 
\end{tikzpicture}\;,
\qquad
\begin{tikzpicture}[scale=0.35,baseline=.6cm]
	\node at (0,0.5)  [root] (root) {};
	\node at (-2,1)  [dot] (left) {};
	\node at (2,1)  [dot] (right) {};
	\node at (-1.5,3)  [dot] (mleft) {};
	\node at (1.5,3)  [dot] (mright) {};
	\node at (0,3) [dot] (top) {};
\node at (0,4.5)  [dot] (ttop) {};
\node at (-2.5,4)  [dot] (ttopleft) {};
\node at (2.5,4)  [dot] (ttopright) {};	
	\draw[kernel]   (root) to (left);
	\draw[kernel]   (root) to (right);
	\draw[kernel]   (left) to (ttopleft);
	\draw[kernel]   (right) to (ttopright);
	\draw[kernel] (right) to  (mright); 
	\draw[kernel] (left) to  (mleft); 
	\draw[arrho] (top) to  (mleft); 
	\draw[arrho] (top) to  (mright); 
	\draw[arrho] (ttop) to  (ttopright); 
	\draw[arrho] (ttop) to (ttopleft);
	\end{tikzpicture} \ .
\end{equation}

Then, one sees exactly\footnote{Here we use that $\rho$ is even.} as in Subsection~\ref{subsec:heat kernel reg} that the only divergent contribution to the integrals in \eqref{eq:general_integral_involving_Gammas}\&\eqref{eq:general_integral_involving_nablaGammas} are given when all occourences of $\Gamma^*$, resp.\ $\nabla\Gamma^*$ are replaced by $Z^*_0$, resp.\ $Z^*_{0,1}$ and that all remaining terms converge as $\epsilon \to 0$. From now on we focus on the
former case, i.e.\ on the integral in \eqref{eq:general_integral_involving_Gammas}, the latter corresponding to \eqref{eq:general_integral_involving_nablaGammas} can be treated analogously.
The counterterm corresponding to the one considered in the previous section is thus given by 
\begin{equ}\label{eq:naive_counterterm}
\mathcal{I}_{\epsilon,\CCG} (z_{\star}) = \int_{z\in (\mathbb{R}\times \mathbb{R}^d)^{\CCV_0}} \prod_{e\in \CCE_+} C(z_\star) w^{(z_{\star})} (z_{e_+},z_{e_-})  \prod_{e\in \CCE_0}  \varrho^\epsilon  (z_{e_+},z_{e_-})  dx\ 
\end{equ}
with $\CCG$ being one of the graphs in \eqref{eq:graphs}.
Furthermore, it is clear subtracting $\mathcal{I}_{\epsilon,\CCG} (z_{\star})$ from \eqref{eq:general_integral_involving_Gammas} corresponds in the limit $\epsilon\to 0$ to the same renormalisation as considered in the previous subsection.

Next we argue that $\mathcal{I}_{\epsilon,\CCG} (z_{\star})$ has the same asymptotic behaviour as the counterterms obtained in Subsection~\ref{subsec:heat kernel reg} (but with the constants $\alpha,\beta_\epsilon$ in general depending additionally on $\rho$). The following consequence of Taylor's theorem will be useful.
\begin{lemma}\label{lem:dif delta approx}
There exists a constant $C$ depending on $\rho$ and $A$, such that
$$|\varrho^\epsilon  (z_1,z_2) - \varrho^{z,\epsilon}  (z_1,z_2) |\leq C \frac{1}{\epsilon^d} \mathbf{1}_{B_\epsilon } (x_1-x_2) \left(|x_1-x| + |t_1-t|^{1/2}\right) \ .$$
\end{lemma}

We define
\begin{equation}\label{eq:counterterm with correct assimptotics}
\mathring{\mathcal{I}}_{\epsilon,\CCG} (z_{\star}):= \int_{(\mathbb{R}\times \mathbb{R}^d)^{\CCV_0}} \prod_{e\in \CCE_+} C(z_\star) w^{(z_{\star})} (z_{e_+},z_{e_-})  \prod_{e\in \CCE_0}  \varrho^{z_{\star}, \epsilon}  (z_{e_+},z_{e_-})  dx \ ,
\end{equation}
and note by a simple substitution that $\mathring{\mathcal{I}}_{\epsilon,\CCG}$ has the same form as counterterms explicitly calculated in Section~\ref{subsec:heat kernel reg}. Therefore, we conclude this section by showing that
\begin{equation}\label{eq:counterterm diff to 0}
| \mathcal{I}_{\epsilon,\CCG} (z_{\star}) - \mathring{\mathcal{I}}_{\epsilon,\CCG} (z_{\star})|\to 0
\end{equation}
(uniformly in $z_\star$) as $\epsilon\to 0$.
Using an analogous decomposition as in \cite[Lemma~A.4]{HQ18} we write
$w^{(z_{\star})} = \sum_{n=0}^\infty w^{(z_{\star})}_n$ and for $\mathbf{n}\in \NN^{\CCE_+}$
$$W^{\mathbf{n}}(z):= \prod_{e\in \CCE_+} C(z_\star) w^{(z_{\star})}_{\mathbf{n}_e} (z_{e_+},z_{e_-})\ , $$
we find that
\begin{align}
\mathcal{I}_{\epsilon,\CCG} (z_{\star}) &- \mathring{\mathcal{I}}_{\epsilon,\CCG} (z_{\star}) \nonumber \\
&= \sum_{\mathbf{n}\in \NN^{\CCE_+}}
\int_{(\mathbb{R}\times \mathbb{R}^d)^{\CCV_0}} W^{\mathbf{n}}(z) \left(\prod_{e\in \CCE_0}  \varrho^{z_{\star}, \epsilon}  (z_{e_+},z_{e_-})  
-\prod_{e\in \CCE_0}  \varrho^{\epsilon}  (z_{e_+},z_{e_-})  \right) \ dz
 \ . \label{eq:diverging counterterm decomp}
\end{align}
To conclude, we apply dominated convergence together with the following observations
\begin{enumerate}
\item Each summand of \eqref{eq:diverging counterterm decomp} converges to $0$ (uniformly in $z_*$).
\item Due to Lemma~\ref{lem:dif delta approx} one finds that for each of the diagrams in \eqref{eq:graphs} for the equations considered, 
\begin{equation}\label{eq:unif bound for dom con}
\sum_{\mathbf{n}\in \NN^{\CCE_+}} \sup_{\epsilon\in (0,1)} \left|
\int_{(\mathbb{R}\times \mathbb{R}^d)^{\CCV_0}} W^{\mathbf{n}}(z) \left(\prod_{e\in \CCE_0}  \varrho^{z_{\star}, \epsilon}  (z_{e_+},z_{e_-})  
-\prod_{e\in \CCE_0}  \varrho^{\epsilon}  (z_{e_+},z_{e_-})  \right) \ dz \right| 
\end{equation}
is bounded uniformly in $\epsilon$ and independently of $z_\star$.
\end{enumerate}

\begin{remark}\label{rem:fail}
Let us emphasise, that whether the uniform bound on \eqref{eq:unif bound for dom con} holds depends crucially on the power counting of the diagram considered. It fails for the ``cherry" of the $\phi^3_4$ equation.
\end{remark}

%
%Crucially, the function $\varrho$ needs to be $\mathfrak{R}_\zeta(x)$- invariant for the analogue of the term of \eqref{eq:vanishing_term_phi^4} to vanish, otherwise additional $A$ dependent counter-terms are needed.
\subsection{Further regularisations}\label{sec:flat_reg}
Finally, we observe how renormalisation functions should be chosen when working with more general regularisations of the noise in order to obtain the same solutions in limit when the regularisation is removed. To this end, let us consider regularised noises $\xi_\eps(z):= \int \rho_\eps(z,z')d\xi(z')$
for some compactly supported continuous function $\rho_\eps:\mathbb{R}^{d+1}\times \mathbb{R}^{d+1}$ converging to the Hausdorff measure on the 
diagonal $\{(z,z')\in \mathbb{R}^{d+1}\times \mathbb{R}^{d+1} : z=z'\}$ (with the obvious adaptation when working with spatial white noise).

Thus, the counterterm to subtract can still be represented by the same directed graph as in the previous section $\CCG= (\CCV, \CCE)$,
but where the occurrences of $\varrho^\epsilon  (z_{e_+},z_{e_-})$ in \eqref{eq:naive_counterterm} are replaced by $\rho_\eps(z_{e_+},z_{e_-})$, namely 
$$\mathcal{I}_{\epsilon,\CCG} (z_{\star}) = \int_{z\in (\mathbb{R}\times \mathbb{R}^d)^{\CCV_0}} \prod_{e\in \CCE_+} C(z_\star) w^{(z_{\star})} (z_{e_+},z_{e_-})  \prod_{e\in \CCE_0}  \rho_\eps(z_{e_+},z_{e_-})  dz. \ 
$$ 
Observe that the dependence of $w^{(z_{\star})}$ on the variable $z_{\star}$ is really only through the coefficient field, i.e.
we can write $w^{(z_{\star})}= \tilde{w}^{(A(z_{\star}))}$ where for $z\neq z'$
the functions $\tilde{w}^{(\cdot)}(z,z')$ map positive definite matrices to the real numbers.
 Thus we find that 
 $$\mathcal{I}_{\epsilon,\CCG} (z_{\star}) = \prod_{e\in \CCE_+} C(z_\star) \cdot \int_{z\in (\mathbb{R}\times \mathbb{R}^d)^{\CCV_0}} \prod_{e\in \CCE_+} \tilde{w}^{A(z_{\star})} (z_{e_+},z_{e_-})  \prod_{e\in \CCE_0}  \rho_\eps(z_{e_+},z_{e_-})  dz. 
$$ 
is of the form $\mathcal{I}_{\epsilon,\CCG} = \tilde{\mathcal{I}}_{\epsilon,\CCG}\circ A$ where the functions $\tilde{\mathcal{I}}_{\epsilon,\CCG}$ depend on the graph and the molifier, but importantly do not depend on the coefficient field. 

\begin{remark}
Let us note that one can represent the functions $\tilde{\mathcal{I}}_{\epsilon,\CCG}$ using the same
diagrams as in \eqref{eq:graphs}, but
with a slight twist in interpretation:  the edges $\tikz[baseline=-0.1cm] \draw[kernel] (0,0) to (1,0);$ now
 represent the function $ \tilde{w}^{(\cdot)}$ which in the upper variable is to be evaluated at the matrix $A(z_\star)$.
\end{remark}

\section{Main results}\label{sec:main results} 

Throughout this section write 
$ L = \sum_{i,j=1}^d a_{i,j}(x,t) \partial_i \partial_j + \sum_{i=1}^d b_i(x,t) \partial_i + c(x,t) \ $
and recall that $A$ denotes the matrix with entries $a_{i,j}$ which (without loss of generality) is assumed to be symmetric. Furthermore, recall that we write $a^{i,j}$ for the entries of $A^{-1}$ and 
$C(x,t)= (4\pi)^{-d/2} \det(A(x,t))^{-1/2}$. 
\begin{theorem}\label{thm:g-pam}
Let $\xi$ be white noise on $\mathbb{T}^2$, and $\xi_\epsilon$ its regularisation as in \eqref{eq:heat kernel regularisation spatial white noise} or 
\eqref{eq:general_regular_spatial} for $\rho$ as therein. Let $f_{ij}, g\in \cC^{\infty}(\RR)$ and $u_0\in \mathcal{C}^\alpha(\mathbb{T}^2)$ for $\alpha>0$.
There exists a sequence $\alpha_\epsilon^{\<Xi2>}, \alpha^{\<b2>}_\epsilon\in \mathbb{R}$ depending on $\rho$ such that for $ u_\epsilon$ satisfying $u_\epsilon(0)=u_0$ and solving 
\begin{equation*}
 (\partial_t- L)u =\sum_{i,j=1}^2 f_{ij}(u) \left( \partial_i u \partial_j u - \alpha^{\<b2>}_\epsilon C(x,t) a^{i,j}(x,t)  g^2(u)\right) + g(u)\left( \xi  -\alpha_\epsilon^{\<Xi2>} C(x,t) g'(u)\right) \ ,
\end{equation*}
there exists a (random) $T>0$ and $u: [0,T]\times \mathbb{T}^2 \to \mathbb{R}$ independent of the choice of $\rho$, such that $u_\epsilon \to u$ uniformly as $\epsilon\to 0$ in probability. 
\end{theorem}
\begin{remark}
One can show that $\alpha_\epsilon^{\<Xi2>}, \alpha^{\<b2>}_\epsilon$ in the theorem above can be written as 
$$\alpha_\epsilon^{\<Xi2>} = 2|\log (\epsilon)|+\beta_{\epsilon}^{\<Xi2>} (\rho) \qquad  \alpha^{\<b2>}_\epsilon = |\log (\epsilon)| +\beta_{\epsilon}^{\<b2>} (\rho) $$
where $\beta_{\epsilon}^{\<Xi2>} (\rho), \beta_{\epsilon}^{\<b2>} (\rho)$ converge as $\epsilon \to 0$, c.f.\ Section~\ref{sec:PAM}.
\end{remark}
\begin{proof}[Proof of Theorem~\ref{thm:g-pam}]
The proof follows along the usual lines when applying the theory of regularity structures \cite{Hai14}. For example one can follow \textit{ad verbatim} the proof of \cite[Theorem~16.1]{HS23m} with the only difference in the argument being the the convergence of the $0$-th Wiener Chaos contributions of $\Pi_x  \<Xi2>$ and $\Pi_x  \<b2>$, which for the heat kernel regularisation \eqref{eq:heat kernel regularisation spatial white noise} is the content of Section~\ref{sec:PAM}, and for more general regularisations is contained in Section~\ref{sec:other regularisation}.
\end{proof}

\begin{theorem}\label{thm:phi4_3}
For $d\in \{2,3\}$,
let $\xi$ denote space-time white noise on $\mathbb{R}\times \mathbb{T}^d$ and $u_0\in \mathcal{C}^{\alpha}(\mathbb{T}^d)$ for $\alpha>-\frac{2}{3}$.  Denote by $\xi_\epsilon$ its regularisation as in \eqref{eq:heat kernel regularisation space time white noise} or \eqref{eq:gen_regularisation} for $\varrho$ as therein. Then, there exist sequences of constants
$$\alpha_{\epsilon}^{\<2>,2}, \alpha_{\epsilon}^{\<2>,3},  \alpha_{\epsilon}^{\<22>}\in \mathbb{R}$$ depending on $\varrho$ such that the following holds.
\begin{itemize}
\item If $d=2$, let $u_\epsilon$ denote the solution to
\begin{equation}\label{phi4}
\partial_t u_\epsilon -L u_\epsilon= -u_\epsilon^3+ 3 \alpha_\epsilon^{\<2>,2} C(x,t) u_\epsilon +\xi_\epsilon\ , \qquad
u_\epsilon(0)=u_0 \ .
\end{equation}
\item If $d=3$, let $u_\epsilon$ be the solution to
\begin{equation}\label{phi4}
\partial_t u_\epsilon -L u_\epsilon= -u_\epsilon^3 + 3 \alpha_\epsilon^{\<2>,3} C(x,t) u_\epsilon - 9 \alpha_{\epsilon}^{\<22>} C(x,t)^2 u_\epsilon +\xi_\epsilon\ , \qquad
u_\epsilon(0)=u_0 \ .
\end{equation}
\end{itemize}
%\begin{equation}\label{phi4}
%\partial_t u_\epsilon -L u_\epsilon= -u_\epsilon |u_\epsilon |^2_h+ 3 C^{(2)}_\epsilon C(x,t) u_\epsilon +\xi_\epsilon\ , \qquad
%u_\epsilon(0)\in \cC^{\alpha} \ ,
%\end{equation}
%respectively,
%\begin{equation}\label{phi4}
%\partial_t u_\epsilon -L u_\epsilon= -u_\epsilon |u_\epsilon |^2_h+ 3 C^{(3)}_\epsilon C(x,t)u_\epsilon - 9 C_\epsilon' C(x,t)^2 u_\epsilon +\xi_\epsilon\ , \qquad
%u_\epsilon(0)\in \cC^{\alpha} \ .
%\end{equation}
In both cases, there exists (random) $T>0$ and $u\in C([0,T], \mathcal{D}'(\mathbb{T}^d))$ such that $u_\epsilon \to u$ as $\epsilon\to 0$ in probability. In particular, the limit is independent of the choice of $\rho$
\end{theorem}

\begin{remark}
We first make two observations particular to the non-translation invariant setting.
\begin{enumerate}
\item The renormalisation functions for $\phi^4$ equations depend only on $|\det (A)|^{-1/2}$, which can be interpreted as a ``volume-term", c.f.\  Remark~\ref{rem:riemannian}.
\item For $d=3$, we in general obtain a \textit{genuine} $2$-dimensional solution family for the $\phi^4_3$ equation, in contrast to the translation invariant setting. One thus expects the analogue statement to also hold for non-translation invariant $\phi^4_3$ measures.
\end{enumerate}
\end{remark}

\begin{remark}
On can actually choose $T=+\infty$ a.s.\ in Theorem~\eqref{thm:phi4_3} by an adaptation of the global existence/coming down from infinity results of \cite{MW17},\cite{MW20},\cite{CMW23}. Furthermore, one can show that 
$$\alpha_\epsilon^{\<2>,2}= |\log(\epsilon)| +\beta_\epsilon^{\<2>,2} (\varrho)$$
and 
$$\alpha_\epsilon^{\<2>,3}= \frac{\alpha(\varrho)}{\epsilon}  +\beta_\epsilon^{\<2>,3} (\varrho), \qquad 
\alpha_{\epsilon}^{\<22>} = c|\log(\epsilon)| +\beta_{\epsilon}^{\<22>} (\varrho),$$
for a universal constant $c\in \mathbb{R}$, some $\varrho$ dependent constant $\alpha(\varrho)\in \mathbb{R}$ and where \linebreak $\beta_\epsilon^{\<2>,2}, \beta_\epsilon^{\<2>,3},\beta_{\epsilon}^{\<22>}\in \mathbb{R}$ converge as $\epsilon\to 0$.
\end{remark}
\begin{proof}[Proof of Theorem~\ref{thm:phi4_3}]
For $d=3$, one can for example follow \textit{ad verbatim} (with some simplifications) the proof of \cite[Theorem16.6]{HS23m}, again with the only difference in the argument being, using the notation therein, the treatment of the function resp.\ the renormalised kernel
$$\begin{tikzpicture}[scale=0.35,baseline=0.3cm]
	\node at (0,0)  [root] (int) {};
	\node at (0,2.5) [bluedot] (top) {};
	\draw[blue, keps] (top) to[bend left=60] (int); 
	\draw[blue, keps] (top) to[bend right=60] (int); 
\end{tikzpicture}\;,\qquad \text{resp.} \qquad 
\begin{tikzpicture}[scale=0.35,baseline=-0.1cm]
	\node at (4,0)  [root] (root) {};
	\node at (0,0)  [root] (middle) {};
	\node at (2,-1.5)  [bluedot] (left) {};
	\node at (2,1.5)  [bluedot] (right) {};
	
	\draw[blue,keps,bend right=30] (left) to  (root);	
	\draw[blue,keps,bend left=30] (left) to  (middle);	
	\draw[blue,keps,bend left=30] (right) to  (root);	
	\draw[blue,keps,bend right=30] (right) to  (middle);	
	\draw[blue,kernel] (middle) to (root);
\end{tikzpicture}\;,
$$ i.e. \cite[Lemma~16.13 \& Lemma~16.15]{HS23m}.
For the heat kernel regularisations \eqref{eq:heat kernel regularisation spatial white noise} the replacement of these Lemmas is the content of Section~\ref{sec:phi43}, and for more general regularisations it is contained in Section~\ref{sec:other regularisation}.

The case $d=2$ follows either as a (simple) special case of argument for $d=3$ or can be obtained by a straightforward adaptation of the argument in \cite{DD03}, see also \cite{TW18}.
\end{proof}

\begin{theorem}\label{thm:KPZ}
Let $\xi$ denote space-time white noise on $\mathbb{T}$ and $\xi_\epsilon$ its regularisation as in \eqref{eq:heat kernel regularisation space time white noise} or \eqref{eq:gen_regularisation} for $\varrho$ as therein.
Let $u\in \cC^\alpha(\mathbb{T})$ for $\alpha>0$, then there exists a sequence of constants 
$\alpha_\epsilon^{\<2>}$, $\alpha_\epsilon^{\<22j>}$, $\alpha_\epsilon^{\<211>}\in \mathbb{R}$
(depending on $\varrho$) such that for $u_\epsilon$ solving 
$$
\partial_t u_\epsilon -L u_\epsilon= (\partial_x u_\epsilon)^2  -\alpha_\epsilon^{\<2>} \cdot a^{-3/2} (x,t) - \left( \alpha_\epsilon^{\<22j>} + 4\alpha_\epsilon^{\<211>} \right) \cdot  a^{-4}  (x,t) +\xi_\epsilon\ , \qquad
u_\epsilon(0)=u_0\ ,
$$
there exists a (random) $T>0$ and $u: [0,T]\times \mathbb{T} \to \mathbb{R}$, such that $u_\epsilon \to u$ uniformly as $\epsilon\to 0$ in probability and, furthermore, u is independent of the choice of $\varrho$.
\end{theorem}

\begin{proof}[Proof of Theorem~\ref{thm:KPZ}]
The argument, for example, adapts \textit{mutatis mutandis} from the proof of \cite[Theorem~16.20]{HS23m}, since the diagrams to renormalise are formally identical to the ones for the $\phi^3_4$ equation therein.  Here the content of Section~\ref{sec:KPZ}, resp. Section~\ref{sec:other regularisation} replaces \cite[Lemma~16.24]{HS23m}. 
\end{proof}

\begin{remark}
Let us emphasise that in the proofs above we could also refer to \cite{Hai14} for the analytic step, to \cite{BB21} for the identification of the renormalised equation (once renormalisation functions are chosen) and to \cite{HS23m} for the construction of the regularity structure, the renormalisation group and the stochastic estimates. For the readers convenience we decided to refer directly to \cite{HS23m} instead, as in particular \cite{BB21} uses somewhat different notation.
\end{remark}
\begin{remark}
One obtains the direct analogues of Theorem~\ref{thm:g-pam}, Theorem~\ref{thm:phi4_3} and Theorem~\ref{thm:KPZ}, if one uses the more general mollifiers 
of Section~\ref{sec:flat_reg} with the only modification that one has to replaces the renormalisation functions in those theorems with the respective functions $ \tilde{\mathcal{I}}_{\epsilon,\CCG}\circ A$
defined in that section.
\end{remark}

\begin{prop}
In the setting of Proposition~\ref{prop continuity} consider the equations of Theorems~\ref{thm:g-pam}, \ref{thm:phi4_3} and \ref{thm:KPZ}. 
Then, for each of these equations if we denote by $u^\lambda$ the solution to the equation with $L$ replaced by $L^\lambda$, the map 
$[0,1]\ni \lambda\mapsto u^\lambda$ is Lipschitz continuous (in the topology stated in the respective theorem and for some (possibly smaller) $T>0$ independent of $\lambda$).
\end{prop}
\begin{proof}
Denote by $Z^\lambda$ the renormalised model associated each $\Gamma^\lambda$.
In view of Theorem~\ref{thm:fixed point} the proof of the proposition amounts to checking that for each respective regularity structure the map
$$[0,1]\ni\lambda\mapsto (K^\lambda, {Z}^\lambda) \in \bfK^\beta_{L,R} \ltimes \mathcal{M}$$ is Lipschitz continuous.
But this, in view of Proposition~\ref{prop continuity}, follows by exactly the same way as when establishing convergence of models when regularising the noise.
\end{proof}
	\bibliographystyle{Martin}
	\bibliography{Homog.bib}

\end{document}